\documentclass[12pt]{article}
\usepackage{amssymb, amsmath,amsthm}
\usepackage[dvips]{graphicx}
\usepackage{color}

\date{}

\begin{document}

\title
{\large{Kato's Inequality for Magnetic Relativistic
Schr\"odinger Operators}\\
\bigskip\bigskip
\normalsize{\it By}\\
\bigskip\bigskip
\normalsize{Fumio Hiroshima$^{1)}$},
\normalsize{ Takashi Ichinose$^{2)}$}
\normalsize{and J\'ozsef L\H{o}rinczi$^{3)}$}\\
\bigskip\bigskip
\small{1) Faculty of Mathematics, Kyushu University, 744 Motooka,
          Fukuoka, 819-0395, Japan; {\tt hiroshima@math.kyushu-u.ac.jp}}\\
\small{2) Department of Mathematics, Kanazawa University,
          Kanazawa, 920-1192, Japan; {\tt ichinose@kenroku.kanazawa-u.ac.jp}}
\\
\small{3) Department of Mathematical Sciences, Loughborough University,
          Loughborough LE11 3TU, United Kingdom; {\tt  J.Lorinczi@lboro.ac.uk}}
\\
}


\maketitle

\abstract
{Kato's inequality is shown for the magnetic relativistic Schr\"odinger
operator $H_{A,m}$ defined as the operator theoretical {\it square root} of
the selfadjoint, magnetic nonrelativistic Schr\"odinger operator
$(-i\nabla-A(x))^2+m^2$ with an $L^{2}_{\text{\rm loc}}$ vector potential $A(x)$. }

\noindent------------------------------



\medskip\noindent
{\bf Running Head}: Kato's Inequality

\bigskip\noindent
{\bf Footnotes for the first page}

\noindent------------------------------

\smallskip\noindent
Communicated by H. Okamoto. Received January 15, 2016. Revised June 13, 2016; July 11, 2016.

\noindent
F. Hiroshima: Faculty of Mathematics, Kyushu University, 744 Motooka,
          Fukuoka, 819-0395, Japan; e-mail: {\tt hiroshima@math.kyushu-u.ac.jp}

\smallskip
\noindent
T. Ichinose: Department of Mathematics, Kanazawa University,
          Kanazawa, 920-1192, Japan; e-mail: {\tt ichinose@kenroku.kanazawa-u.ac.jp}

\smallskip
\noindent
J. L\H{o}rinczi: Department of Mathematical Sciences, Loughborough University,
          Loughborough LE11 3TU, United Kingdom; e-mail: {\tt  J.Lorinczi@lboro.ac.uk}

\bigskip\bigskip\bigskip\noindent

\noindent+++++++++++++++++++++++++++++++++++++

\noindent
The Corresponding Author:

Takashi ICHINOSE,$\,\,\,$ E-mail: ichinose@kenroku.kanazawa-u.ac.jp,

Tel: +81-76-244-8631; $\,\,\,$ Fax: +81-76-264-5738

\makeatletter
\renewcommand{\theequation}{%
  \thesection.\arabic{equation}}
\@addtoreset{equation}{section}
\makeatother

\newtheorem{dfn}{Definition}[section]
\newtheorem{thm}[dfn]{Theorem}
\newtheorem{prp}[dfn]{Proposition}
\newtheorem{lemma}[dfn]{Lemma}

\newpage
\section{Introduction} 
Consider the magnetic relativistic Schr\"odinger operator
\begin{equation}
H_{A,m} := \sqrt{(-i\nabla-A(x))^2+m^2}
\end{equation}
in $d$-dimensional space ${\Bbb R}^d$ with vector potential
$A(x) := (A_1(x), \dots, A_d(x))$ and rest mass $m\geq 0$,
which may be thought of being a quantum Hamiltonian
corresponding to the classical relativistic Hamiltonian symbol
$\sqrt{(\xi-A(x))^2 + m^2}, \,\, (\xi,x) \in {\Bbb  R}^d \times {\Bbb  R}^d$.
It is known that when $A(x)$ is an ${\Bbb R}^d$-valued function
belonging to
{$[L^2_{\text{\rm loc}}({\Bbb R}^d)]^d
\equiv L_{\text{\rm loc}}^2({\Bbb  R}^d; {\Bbb  R}^d)$,
}
it becomes a selfadjoint operator in $L^2({\Bbb R}^d)$, which is essentially
selfadjoint on $C_0^{\infty}({\Bbb R}^d)$ so that $H_{A,m}$ has a domain
containing $C_0^{\infty}({\Bbb R}^d)$ as an operator core (e.g see [CFKiSi87, p.9]).
We shall assume that $d\geq 2$, since in case $d=1$ any magnetic vector potential
can be removed by a gauge tranformation.
For $A=0$ we put $H_{0,m} = \sqrt{-\Delta +m^2}$, where $-\Delta$ is
the {\it minus-signed} Laplacian
{$-\big(\frac{\partial^2}{\partial x_1^2} +\cdots +
 \frac{\partial^2}{\partial x_d^2}\big)$
} as well as
a nonnegative selfadjoint operator realized in $L^2({\Bbb R}^d)$
having the Sobolev space $H^2({\Bbb R}^d)$ as its domain.

The aim of this paper is to show Kato's inequality for this magnetic
relativistic Schr\"odinger operator $H_{A,m}$ or $H_{A,m}-m$,
when $A$ is an ${\Bbb R}^d$-valued $L^{2}_{\text{\rm loc}}$ function
in ${\Bbb R}^d$.

\begin{thm} (Kato's inequality).
Let $m\geq 0$ and assume $A \in [L^{2}_{\text{\rm loc}}({\Bbb R}^d)]^d$.
If $u \in L^2({\Bbb R}^d)$ with
$H_{A,m} u \in L^1_{\text{\rm loc}}({\Bbb R}^d)$,
then the following distributional inequality holds:
\begin{align}%
 &\hbox{\rm Re}[(\hbox{\rm sgn}\, u) H_{A,m} u] \geq H_{0,m} |u|,\\
\intertext{\rm or}
 &\hbox{\rm Re}[(\hbox{\rm sgn}\, u) [H_{A,m}-m] u] \geq [H_{0,m}-m]|u|.
\end{align}
Here $\hbox{\rm sgn}$ is a bounded function in ${\Bbb R}^d$ defined by
$$
 (\hbox{\rm sgn}\,u)(x) =
 \left\{
 \begin{array}{ll}
  \overline{u(x)}/|u(x)|, &\hbox{\rm if}\,\, u(x) \not= 0,\\
       0,                 &\hbox{\rm if}\,\, u(x) = 0.
 \end{array}\right.
$$
\end{thm}

\bigskip
Note here that $H_{A,m} u$ with $u \in L^2({\Bbb R}^d)$ makes sense as
a distribution in ${\Bbb R}^d$ (for this, see Lemma 2.2 with $\alpha=1$
and a few lines after its proof).
{ A characteristic feature in this situation is
} that $H_{A,m}$ is a nonlocal operator
defined by the operator-theoretical square root of a nonnegative
selfadjoint operator. It is not a differential operator, and neither
an integral operator nor a pseudo-differential operator associated
with a certain tractable symbol.
The point which becomes crucial is in how to go without knowledge on
regularity of the weak solution $u \in L^2({\Bbb R}^d)$ of equation
$H_{A,m} u = f$ for a given $f \in L^1_{\text{\rm loc}}({\Bbb R}^d)$.
Thus the present inequality (1.2)/(1.3) differs from an abstract form of 
Kato's inequality such as in [Si77]
{by being substatially sharp. }

\bigskip
An immediate corollary is the following theorem, which has been known
(e.g. [FLSei08], [HILo12]; cf. [I93]).

\begin{thm} (Diamagnetic inequality) 
Let $m\geq 0$ and assume that $A \in [L^{2}_{\text{\rm loc}}({\Bbb R}^d)]^d$.
Then it holds that for $f,\, g \in L^2({\Bbb R}^d)$,
\begin{equation}
|(f, e^{-t[H_{A,m}-m]}g)| \leq (|f|, e^{-t[H_{0,m}-m]}|g|).
\end{equation}
\end{thm}

\bigskip
Once Theorem 1.1 is established, we can apply it to show the following theorem
on essential selfadjointness of the relativistic Schr\"odinger operator
with both vector and scalar potentials $A(x)$ and $V(x)$:
\begin{equation}%
 H_{A,V,m} := H_{A,m} + V.
\end{equation}

\begin{thm} 
Let $m\geq 0$, assume that $A \in [L_{\text{\rm loc}}^{2} ({\Bbb R}^d)]^d$ and
let $V \in L_{\text{\rm loc}}^{2} ({\Bbb R}^d)$ with
$V(x) \geq 0$ a.e. Then $ H_{A,V,m} = H_{A,m} + V$ is essentially selfadjoint
on $C_0^{\infty}({\Bbb R}^d)$ and its unique selfadjoint extension is bounded
below by $m$.
\end{thm}

\bigskip
We shall show inequality (1.2)/(1.3), basically along the idea and
method of Kato's original proof in [K72] for the magnetic
{\it nonrelativistic} Schr\"odinger operator $\frac12(-i\nabla-A(x))^2$.
As a matter of fact, we follow  the method of proof
modified for the existing form of Kato's inequality in [I89], [ITs92]
for {\it another} magnetic relativistic Schr\"odinger operator which is
defined as a Weyl pseudo-differential operator associated with the same
relativistic classical symbol $\sqrt{(\xi-A(x))^2 + m^2}$.
{ However, this is not sufficient, and we need further modifications
{\it using operator theory},
since pseudo-differential calculus does not seem useful.
Starting from the assumption of the theorem that $u \in L^2$ and
$H_{A,m}u \in L^1_{\text{\rm loc}}$, it appears to be impossible to show
the regularity of $u$ that
{ $\partial_j u \in L^1_{\text{\rm loc}}\,, \, 1\leq j \leq d$, and/or
$H_{0,m}u \in L^1_{\text{\rm loc}}\,$,
} which may be due to the fact that the operators
$\partial_j\cdot (-\Delta+m^2)^{-1/2}$, $1\leq j \leq d$,
 are {\it not} bounded from $L^1$ to $L^1$,
though they are {\it bounded} from $L^1$ to {\it weak $L^1$-space}.
Therefore we make a detour by going via the case of the fractional power
$(H_{A,m})^{\alpha}$ with $\alpha<1$
.
Verifying that the assumption implies that
$(H_{A,m})^{\alpha}u \in L^1_{\text{\rm loc}}$ for $0<\alpha<1$,
 we show the asserted inequality
{\it first for the case $0<\alpha<1$}, i.e. inequality (1.2)/(1.3)
with the pair $H_{A,m}\,$, $H_{0,m}$, replaced
by the pair $(H_{A,m})^{\alpha}\,$, $(H_{0,m})^{\alpha}$, respectively, and
{\it then for the case $\alpha=1$}, appealing to the fact, to be shown,
that $(H_{A,m})^{\alpha}u$ converges to
$H_{A,m}u$ in $L^1_{\text{\rm loc}}$ as $\alpha \uparrow 1$
.
The proof is presented separately according to $m>0$ and $m=0$,
in a self-contained manner.
}

A comment is in order on our starting assumption for $u$,
namely, why the theorem is formulated with assumption
that $u \in L^2$ and $H_{A,m}u \in L^1_{\text{\rm loc}}$,
but not that {\it both $u$ and $H_{A,m}u$ are $L^1_{\text{\rm loc}}$}.
For this question, recall that the original form of Kato's inequality 
for nonrelativistic Schr\"odinger operators $\frac12(-i\nabla-A(x))^2$
is formulated under the  assumption
that {\it both $u$ and $\frac12(-i\nabla-A(x))^2 u$ are
$L^1_{\text{\rm loc}}$}.
The answer is simply because of avoiding inessential complexity
coming from the fact that  $H_{A,m}$ is a nonlocal operator.

The relativistic Schr\"odinger operator $H_{0,m} = \sqrt{-\Delta+m^2}$
without vector potential was first considered in [W74] and [He77] for
spectral problems. The magnetic relativistic Schr\"odinger operator
$H_{A,m}$ like (1.1) is used to study problems related to 
``stability of matter" in
relativistic quantum mechanics in [LSei10].
On the other hand, a problem of representing by path integral
the relativistic Schr\"odinger semigroup with generator $H_{A,m}$ has been
also studied. A result is a formula of Feynman--Kac--It\^o type
(cf. [Si79/05]), earlier in [DeRiSe91], [DeSe90] and also in [N00],
which has been recently extensively studied in [HILo12], [HILo13] (cf. [LoHBe11]).
The problem is connected with a L\'evy process obtained by subordinating Brownian
motion ([Sa99], [Ap04/09]). A weaker version of Kato's inequality as well
as the diamagnetic inequality was given in our paper [HILo12], to which
the present work adds further results.

In Section 2 some technical lemmas are given,
which are used in the proof of theorems.
They concern some basic inequalities in $L^2$ and $L^p$
connected with the semigroups and/or inverse (resolvent)
for the magnetic {\it nonrelativistic} (but not {\it relativistic})
Schr\"odinger operator $(-i\nabla-A)^2 +m^2$,
which is the square of our magnetic relativistic Schr\"odinger operator
$H_{A,m}$.
For the sake of regularization of $H_{A,m}$, its fractional
powers $(H_{A,m})^{\alpha}$ with $0<\alpha<1$ are also considered
through the semigroup of the magnetic {\it nonrelativistic} Schr\"odinger
operator to estimate in local $L^1$-norm a kind of
difference, being a distance in a particular sense, 
between $(H_{A,m})^{\alpha}$ and $(H_{0,m})^{\alpha}$,
 each applied to a function.

In Section 3 we prove the theorems.
 Section 4 is to make concluding remarks about how 
 the issue is going with the other two magnetic relativistic Schr\"odinger 
operators associated with the same symbol.
Appendix A provides for an explicit expression of the integral kernel
(heat kernel) of the semigroup $e^{-t[(H_{0,m})^{\alpha}-m^{\alpha}]}$
for the free fractional power $(H_{0,m})^{\alpha}$ together with
the density (function) of the associated L\'evy measure $n^{m,\alpha}(dy)$.
For basic facts on the magnetic relativistic Schr\"odinger operator, we
 refer, e.g., to [LLos01] and [BE11].

Finally, we note that we have defined the fractional powers of
$H_{A,m}$ mainly through the magnetic nonrelativistic
Schr\"odinger semigroup. However, an alternative way is
to define them through the Dunford integral via the
{\it resolvent of the magnetic nonrelativistic
Schr\"odinger operator}.

\section{Technical Lemmas} 
Throughout this paper, we denote by $(\cdot, \cdot)$ the Hilbert space
inner product which is
{ sesquilinear, i.e.}
conjugate-linear in the first argument and linear in the second
(the physicist's convention),
and by $\langle\cdot, \cdot\rangle$
the bilinear inner product which is linear in both the arguments.

Our main object of study is the operator 
$H_{A,m} := [(-i\nabla -A)^2+m^2]^{\frac12}$
in (1.1) with assumption that $A \in [L^2_{\text{\rm loc}}({\Bbb R}^d)]^d$,
which is a selfadjoint operator in $L^2({\Bbb R}^d)$ defined as the square root
of the nonnegative selfadjoint (Schr\"odinger) operator
$(-i\nabla -A)^2+m^2$ in $L^2({\Bbb R}^d)$.
For $m=0$, $H_{A,0} = |-i\nabla -A|$.
Among them, the following identity holds:
\begin{eqnarray}%
 \|H_{A,m} u\|_{L^2}^2
&\!\!=\!& \big( u, (H_{A,m})^{2}u \big)
 = \big( u, [(-i\nabla -A)^2 + m^2]u \big) \nonumber\\
&\!\!=\!& \sum_{j=1}^d \|(-i\partial_j -A_j)u\|_{L^2}^2 + m^2\|u\|_{L^2}^2
  = \|H_{A,0} u\|_{L^2}^2 + m^2\|u\|_{L^2}^2\,,   
\end{eqnarray}
with $u \in C_0^{\infty}({\Bbb R}^d)$ for all the five members and
with $u$ in the domain of $H_{A,m}$ for the first, fourth and fifth
 members. The nonrelativistic Schr\"odinger op
erator $(-i\nabla -A)^2+m^2$ concerned
is the selfadjoint operator associated with this quadratic form (2.1), which
 has $C_0^{\infty}({\Bbb R}^d)$ as a form core (e.g [CFKiSi87, 1.3, pp.8--9]).
As a result, $H_{A,m}$ has
$C_0^{\infty}({\Bbb R}^d)$ as an operator core, in other words,
$H_{A,m}$ is a nonnegative selfadjoint operator in $L^2({\Bbb R}^d)$
having domain $D[H_{A,m}]
:= \{u \in L^2({\Bbb R}^d);\; (i\partial_j +A_j)u \in L^2({\Bbb R}^d),
\,\, \partial_j := \partial/\partial_{x_j},\,\,  1\leq j \leq d \}$,
being essentially selfadjoint on $C_0^{\infty}({\Bbb R}^d)$.
Though $i\nabla+A \equiv (i\partial_1+A_1, \dots, i\partial_d+A_d)$
is a closed linear operator of $[L^2({\Bbb R}^d)]^d$ into itself with domain
$D[i\nabla + A] :=
\{(u_1, \dots, u_d) \in [L^2({\Bbb R}^d)]^d;\;
  (i\partial_j +A_j)u \in L^2({\Bbb R}^d),
\,\, \partial_j := \partial/\partial_{x_j},\,\,  1\leq j \leq d \}$,
we will also abuse notation to write {\sl the first term of
the fourth member} of (2.1) as $\|(-i\nabla-A)u\|_{L^2}^2$.

For the proof of Theorem 1.1, however, we need to consider $H_{A,m}$ also
on $L^p$ spaces, and further the fractional powers
$(H_{A,m})^{\alpha},\, 0<\alpha<1$ of $H_{A,m}$.
The aim of this section concerns the issue such as some estimates
connected with them.

{\it As for the constant $m$, unless otherwise stated,
we assume in this section that $m>0$, and keep
{ assuming}
it also in Section 3, until we come to consider the case including $m=0$
 at the final stage of the proof of Theorem 1.1. } Therefore, in case $m>0$,
$H_{A,m}$  has bounded inverse $(H_{A,m})^{-1}$,  as well as
$[(-i\nabla -A)^2 + m^2]$ has bounded inverse  $[(-i\nabla -A)^2 + m^2]^{-1}$.

\subsection{Some inequalities related to magnetic
nonrelativistic Schr\"odinger operators on $L^p$}

The operators $H_{A,m}$ may be considered not only in $L^2$ but also
in $L^p,\,\,1\leq p<\infty$, in particular, for $p=1$.
The square of $H_{A,m}$ becomes a magnetic nonrelativistic Schr\"odinger
operator $(-i\nabla-A)^2 +m^2$.
Some basic inequalities are given which are related to the
 magnetic nonrelativistic Schr\"odinger semigroup
$e^{-t(H_{A,m})^2}$ and inverse (resolvent) $((H_{A,m})^2)^{-1}$
on $L^p$, though {\it not} with the magnetic relativistic Schr\"odinger
semigroup $e^{-tH_{A,m}}$ and inverse (resolvent) $(H_{A,m})^{-1}$.
They will be useful throughout the paper.

In the beginning, let us repetitively confirm the notations to be used:
\begin{align*}
(H_{A,m})^2 &= (-i\nabla-A)^2 +m^2,\quad
(H_{A,0})^2 = (-i\nabla-A)^2 ,\\
(H_{0,m})^2 &=  -\Delta +m^2, \quad (H_{0,0})^2 = -\Delta.
\end{align*}

\begin{lemma}
Suppose $A \in [L^2_{\text{\rm \rm loc}}({\Bbb R}^d)]^d$.
Then the following inequalities hold.

(i) Let $1\leq p \leq \infty$. For $m\geq 0$,
\begin{align*}%
\|e^{-t(H_{A,m})^2}\|_{L^p \rightarrow L^p}
  &\leq \|e^{-t(H_{0,m})^2}\|_{L^p \rightarrow L^p}
  \equiv  \|e^{-t(-\Delta+m^2)}\|_{L^p \rightarrow L^p}
  \leq e^{-m^2 t} \leq 1, \quad t>0.
\end{align*}%
For $m>0$ and $\beta>0$,
\begin{align*}%
\|((H_{A,m})^2)^{-\beta}\|_{L^p \rightarrow L^p}
 &\leq \|((H_{0,m})^2)^{-\beta}\|_{L^p \rightarrow L^p}, \qquad t>0.
\end{align*}%

(ii) Let $1\leq p < \infty$. The operators 
$e^{-t(H_{0,0})^2} (-i\nabla)$ and $e^{-t(H_{0,0})^2} (-\Delta)$
can be extended to be bounded operators
on $[L^p({\Bbb R}^d)]^d$ and $L^p({\Bbb R}^d)$:
\begin{align*}%
\|e^{-t(H_{0,0})^2} (-i\nabla)\|_{[L^p]^d \rightarrow [L^p]^d}
       \leq C_{1p} t^{-1/2},
\quad
\|e^{-t(H_{0,0})^2} (-\Delta)\|_{L^p \rightarrow L^p}
       \leq C_{2p} t^{-1},  \quad t>0,
\end{align*}%
with constants $C_{1p} >0$ and $C_{2p}$ independent of $t$.

\medskip
(iii)  Let $m\geq 0$. The operators  $H_{A,m} e^{-(H_{A,m})^2}$ and
$(H_{A,m})^2 e^{-t(H_{A,m})^2}$ can be extended to be bounded operators
on $L^2({\Bbb R}^d)$: 
$$
\|H_{A,m} e^{-t(H_{A,m})^2}\|_{L^2 \rightarrow L^2} \leq (2et)^{-1/2},
\quad \|(H_{A,m})^2 e^{-t(H_{A,m})^2}\|_{L^2 \rightarrow L^2} 
                                 \leq (et)^{-1},  \quad t>0.
$$

(iv) The operators $e^{-t(-i\nabla -A)^2}(i\nabla + A)$ and
$(i\nabla + A)e^{-t(-i\nabla -A)^2}$
can be extended to be bounded operators on $[L^2({\Bbb R}^d)]^d$: 
\begin{align*}%
&\|e^{-t(-i\nabla -A)^2}(i\nabla + A)\|_{[L^2]^d \rightarrow [L^2]^d}
      \leq \big(\tfrac{d}{2et}\big)^{\frac12},\\
&\|(i\nabla + A)e^{-t(-i\nabla -A)^2}\|_{[L^2]^d \rightarrow [L^2]^d}
      \leq \big(\tfrac{d}{2et}\big)^{\frac12}, \quad t>0.
\end{align*}%
\end{lemma}

The assertion (ii) of Lemma 2.1 may be an $L^p$ version of (iii) or (iv) 
above, though only for a special case of the minus-signed Laplacian 
$-\Delta$ without vector potential $A(x)$.

\bigskip
{\it Proof of Lemma 2.1}. (i) 
This is due to the ingenious
observation given for the magnetic nonrelativistic Schr\"odinger operator
$(-i\nabla-A(x))^2$ with
{ $A \in L^2_{\text{\rm loc}}$
}
in [Si79, Theorem 2.3, p.40], [Si82, Sect. B13, p.490], since
$(H_{A,m})^2 = (-i\nabla-A(x))^2+m^2$ is nothing but a magnetic (nonrelativistic)
Schr\"odinger operator plus the constant $m^2$.
Following the arguments there we have,
for $1\leq p< \infty$ and for every $u\in C_0^{\infty}({\Bbb R}^d)$,
\begin{align*}
&|e^{-t(H_{A,m})^2}u| \leq e^{-t(H_{0,m})^2}|u| = e^{-m^2t}e^{-t(-\Delta)}|u|,
   \quad \text{\rm  pointwise\,\, a.e.},\\
\intertext{so that
$e^{-(H_{A,m})^2}L^p({\Bbb R}^d)
                 \subseteq L^{\infty}({\Bbb R}^d)\cap L^p({\Bbb R}^d)$,
in fact, for $u \in L^p({\Bbb R}^d)$,}
&\|e^{-t(H_{A,m})^2}u\|_{L^p} \leq e^{-m^2t}\|e^{-t(-\Delta)}|u|\|_{L^p}
  \leq e^{-m^2t}\|u\|_{L^p} \leq \|u\|_{L^p}, \quad t\geq 0.
\end{align*}
Thus we can consider $e^{-t(H_{A,m})^2}$ also as a bounded linear operator mapping
$L^p({\Bbb R}^d)$ into itself. Moreover, it is seen it is a contraction semigroup.
We may use the notations $(H_{A,m})^2, \,\, H_{A,m}$ also to mean
operators $(H_{A,m})^2_p,\,\, (H_{A,m})_p$ in $L^p$ when there is no risk of confusion. Furthermore, 
for the crucial assertion (i), we refer to [Si82, Corollary B.13.3, p.491].

\medskip
(ii) 
In fact, $e^{-t(-\Delta)}$ becomes a holomorphic semigroup on
$L^p({\Bbb R}^d),\,\, 1\leq p <\infty$, for $\text{\rm Re}\, t>0$.
Then for any $f \in L^p({\Bbb R}^d)$,
$v(t) := e^{-t(-\Delta)}f$ gives a unique solution of the heat equation
$\frac{\partial}{\partial t} v(t) = \Delta v(t)$
(see e.g. [K76, IX.\S1.8, p.495] and [K76, IX.\S1.6, Remark 1.22, p.492]).
This implies that $e^{-t(-\Delta)}$ has range in the domain
$D[(-\Delta)]$ of $(-\Delta)$, equivalently, that
$te^{-t(-\Delta)}(-\Delta)$ is uniformly bounded
from $L^p({\Bbb R}^d)$ into itself for every real $t>0$, and
so is $t^{1/2}e^{-t(-\Delta)}(-i \partial_j)$ for each $j = 1,2,\dots, d$.

\medskip
(iii) 
For functions in in $L^2$, the assertion are evident by the spectral theorem,
because $(H_{A,m})^2$ and $H_{A,m}$ are nonnegative selfadjoint operators in
the Hilbert space $L^{2}({\Bbb R}^d)$. Indeed, it is easy to see that
 for $u \in C_0^{\infty}({\Bbb R}^d)$,
\begin{align*}%
&\|e^{-t(H_{A,m})^2}H_{A,m}u\|_{L^2}^2 = (u,(H_{A,m})^2 e^{-2t(H_{A,m})^2} u)
 \leq \sup_{\lambda \geq 0} \lambda e^{-2t\lambda}\|u\|_{L^2}^2
  = (2et)^{-1}\|u\|_{L^2}^2,\\
&\|e^{-t(H_{A,m})^2} (H_{A,m})^2u\|_{L^2}^2 = (u, (H_{A,m})^4 e^{-2t(H_{A,m})^2} u)
 \leq \sup_{\lambda \geq 0} \lambda^2 e^{-2t\lambda}\|u\|_{L^2}^2
  = (et)^{-2}\|u\|_{L^2}^2.
\end{align*}%
This shows (iii).

\medskip
(iv) ({\it proof})
These inequalities follow from (ii). Indeed, for the first one, since
$$
\|e^{- t(-i\nabla -A)^2}(i\nabla +A)\varphi\|_{L^2}^2
 := \sum_{j=1}^d \|e^{- t\sum_{k=1}^d(-i\partial_k-A_k)^2}(i\partial_j + A_j)
   \varphi_j \|_{L^2}^2
$$
for $\varphi = (\varphi_1,\dots, \varphi_d)
\in [C_0^{\infty}({\Bbb R}^d)]^d$, we have only to show that for each $j$
$$
\|e^{- t(-i\nabla-A)^2}(i\partial_j + A_j)\varphi_j \|_{L^2}^2
   \leq (2et)^{-1}\|\varphi_j \|_{L^2}^2.
$$
This is seen as follows: For $m>0$, we have by (ii)
\begin{align*}
&\|e^{- t(-i\nabla-A)^2}(i\partial_j + A_j)\varphi_j \|_{L^2}^2\\
&= e^{2m^2}\|[e^{- tH_{A,m}^2}H_{A,m}][H_{A,m}^{-1}
  ((i\partial_j+A_j)^2+m^2)^{1/2}][((i\partial_j+A_j)^2+m^2)^{-1/2}
  (i\partial_j + A_j)]\varphi_j \|_{L^2}^2\\
&\leq e^{2m^2}(2et)^{-1}\|[H_{A,m}^{-1}
  ((i\partial_j+A_j)^2+m^2)^{1/2}] [((i\partial_j+A_j)^2+m^2)^{-1/2}
  (i\partial_j + A_j)]\varphi_j \|_{L^2}^2\\
&\leq e^{2m^2}(2et)^{-1}\|\varphi_j \|_{L^2}^2.
\end{align*}
Letting $m\downarrow 0$, we have the result.

The second one is shown similarly. This shows (iv),
ending the proof of Lemma 2.1.
\qed

\bigskip\noindent
{\it Remark}.  Nontriviality of the assertion (ii) of this lemma lies in that
$i\nabla+A$ {\it does not commute} with the operator
$(i\nabla+A(x))^2 = \sum_{j=1}^d (i\partial_j+A_j(x))^2$ or $(H_{A,m})^2$

\subsection
{Estimate of a kind of difference between $(H_{A,m})^{\alpha}$ and
$(H_{0,m})^{\alpha}$ in local $L^1$-norm}

In this subsection, we consider the operators  given by the fractional powers
$(H_{A,m})^{\alpha} := [(-i\nabla -A)^2+m^2]^{\alpha/2}, \,\, 0<\alpha \leq 1$,
and provide several lemmas to estimate in local $L^1$-norm a kind of difference
between $(H_{A,m})^{\alpha}$ and $(H_{0,m})^{\alpha}$, each applied to
a function $u$. They are needed to prove Theorem 1.1.
Of course, the case for $\alpha=1$ turns out to be our
operator itself:
$(H_{A,m})^1 \equiv H_{A,m} = [(-i\nabla -A)^2+m^2]^{\frac12}$.

\medskip
Given a positive self-adjoint
operator $S$ in a Hilbert space $L^2({\Bbb R}^d)$ with domain $D[S]$,
we adopt the following definition of its fractional powers $S^{\alpha}$
to be suggested from the identity for the gamma function $\Gamma(\beta)$,
$s^{-\beta} = \frac1{\Gamma(\beta)} \int_0^{\infty} t^{\beta-1}e^{-st}\, dt$
with $t>0$ and $0<\beta \leq 1$ :
for $0 \leq\alpha<1$,
$$
 S^{\alpha}u
  = S^{-(1-\alpha)}\cdot S u
  =\frac1{\Gamma(1-\alpha)} \int_0^{\infty} t^{-\alpha} e^{-tS}\, Su\, dt\,,
  \qquad u \in D[S].
$$

We shall use these formulas, taking for $S$  the nonrelativistic
Schr\"odinger operator
$[(-i\nabla -A)^2 +m^2] = (H_{A,m})^2$ and/or $[-\Delta +m^2] = (H_{0,m})^2$,
but not the relativistic Schr\"odinger operator $H_{A,m}$ and/or $H_{0,m}$.
Thus for $f \in L^2({\Bbb R}^d)$,
\begin{eqnarray}
(H_{A,m})^{-\beta}f
 &=& [(-i\nabla -A)^2+m^2]^{-\frac{\beta}2}f \nonumber\\
 &=& \frac1{\Gamma(\frac{\beta}2)}\int_0^{\infty}
       t^{\frac{\beta}2-1}e^{-t[(-i\nabla -A)^2+m^2]} f dt
       \quad  (0<\beta\leq 2)\,,
\end{eqnarray}%
and similarly for $(H_{0,m})^{-\beta} \equiv [-\Delta+m^2]^{-\beta/2}$
in case $A=0$.
Therefore, for $u \in C_0^{\infty}({\Bbb R}^d)$,  we have
\begin{eqnarray}%
(H_{A,m})^{\alpha}u
&=& [(-i\nabla-A)^2+m^2]^{\alpha/2}u \nonumber\\
&=& \frac1{\Gamma(\frac{2-\alpha}2)}\int_0^{\infty}
    t^{\frac{2-\alpha}2-1}e^{-t[(-i\nabla-A)^2+m^2]}[(-i\nabla-A)^2+m^2]u\, dt
         \nonumber\\
&=& \frac1{\Gamma(\frac{2-\alpha}2)}\int_0^{\infty}
                 t^{-\frac{\alpha}2}e^{-t(H_{A,m})^2}(H_{A,m})^2u\, dt,
    \quad (0 \leq \alpha < 2)\,,
\end{eqnarray}
for $u$ in the domain of $(H_{A,m})^2$,
and similarly for $(H_{0,m})^{\alpha} \equiv [-\Delta+m^2]^{\alpha/2}$
in case $A=0$.
Here note that $H_{A,m}$/$H_{0,m}$, as well as
$S = (-i\nabla -A)^2 +m^2$/ $(-\Delta+m^2)$, has bounded inverse, since
we are assuming in this section that $m>0$.
{ It may be instructive to recognize that for $0<\alpha <1$
the last integral of (2.3) exists
not only for $u \in D[(H_{A,m})^2]$ but also for $u \in D[H_{A,m}]$,
because by Lemma 2.1(iii)
$$
t^{-\frac{\alpha}2}\|e^{-t(H_{A,m})^2}(H_{A,m})^2u\|_{L^2}
\leq t^{-\frac{\alpha}2}\|e^{-t(H_{A,m})^2}H_{A,m}\|\,\|H_{A,m}u\|_{L^2}
= O(t^{-\frac{(1+\alpha)}2}).
$$
}

\begin{lemma}
Let $0<\alpha \leq 1$.
 Assume that $A \in [L^2_{\text{\rm loc}}({\Bbb R}^d)]^d$.
Then if $\varphi \in C_0^{\infty}({\Bbb R}^d)$, then
$(H_{A,m})^{\alpha} \varphi \in L^2({\Bbb R}^d)$.
In fact, it holds for every compact subset $K$ in ${\Bbb R}^d$ that
{\begin{equation}
 \|(H_{A,m})^{\alpha} \varphi\|_{L^2}
\leq [|K|^{\frac12}\big[((m^2+1)^{\frac12}+1)
       +\|\,|A|\,\|_{L^2(K)}]\,
      \big[\|\nabla \varphi\|_{L^{\infty}(K)}+\|\varphi\|_{L^{\infty}(K)}\big],
\end{equation}%
}
for all $\varphi \in C_0^{\infty}({\Bbb R}^d)$ with
$\text{supp}\,\varphi \subseteq K$,
where $|K|$ denotes the volume (Lebesgue measure) of $K$.
\end{lemma}

\bigskip
{\it Proof}.
 Let $\varphi \in C_0^{\infty}({\Bbb R}^d)$ with
$\text{\rm supp}\, \varphi \subseteq K$. Then for $0<\alpha \leq 1$, we have
\begin{eqnarray}
\|(H_{A,m})^{\alpha} \varphi\|_{L^2}^2
&=& \big( \varphi, (H_{A,m})^{2\alpha}\varphi \big)
= \big( \varphi, [(-i\nabla -A)^2 + m^2]^{\alpha}\varphi \big)\nonumber\\
&\leq& \big( \varphi, [(-i\nabla -A)^2 + m^2+1]^{\alpha}\varphi \big)
                                               \nonumber\\
&\leq& \big( \varphi, [(-i\nabla -A)^2 + m^2+1]\varphi \big)\nonumber\\
&=& \|(-i\nabla -A)\varphi\|_{L^2}^2 + (m^2+1)\|\varphi\|_{L^2}^2
 = \|H_{A,(m^2+1)^{\frac12}}\varphi\|_{L^2}^2\,.
\end{eqnarray}
Here for the first term of the last second member recall our
informal notation mentioned after (2.1).
 Hence
\begin{align*}%
&\|(H_{A,m})^{\alpha} \varphi\|_{L^2}\\
&\leq \|\nabla \varphi\|_{L^2} + \|A \varphi\|_{L^2}
                   + (m^2+1)^{\frac12}\|\varphi\|_{L^2}\\
&\leq |K|^{1/2}\|\nabla \varphi\|_{L^{\infty}(K)}
      + \|\,|A|\,\|_{L^2(K)}\|\varphi\|_{L^{\infty}(K)}
       + (m^2+1)^{\frac12}|K|^{1/2}\|\varphi\|_{L^{\infty}(K)} < \infty,
\end{align*}%
which is finite by assumption on $A$ and $\varphi$.
This shows the desired assertion. \qed

\bigskip
By this lemma, for $0<\alpha \leq 1$
we can {\it define} a distribution $(H_{A,m})^{\alpha} u$ for
$u \in L^2({\Bbb R}^d)$ by
{
$$\langle\,(H_{A,m})^{\alpha} u, \phi\,\rangle
= \langle\,u, (H_{-A,m})^{\alpha} \phi\,\rangle
 = \int (u(H_{-A,m})^{\alpha}\phi )(x) dx,
$$
or
$$(\,(H_{A,m})^{\alpha} u, \phi\,)
=  (\,u, (H_{A,m})^{\alpha} \phi\,)
 =\int (\bar{u}(H_{A,m})^{\alpha}\phi )(x) dx,
$$
for $\phi \in C_0^{\infty}({\Bbb R}^d)$, because, for every compact set $K$
in ${\Bbb R}^d$, we have
\begin{eqnarray*}%
|(\,(H_{A,m})^{\alpha} u, \phi \,)|
&=& |(\, u, (H_{A,m})^{\alpha} \phi \,)|
 \leq \|u\|_{L^2} \|(H_{A,m})^{\alpha} \phi \|_{L^2}\\
&\leq& \|u\|_{L^2} \big[(|K|^{\frac12}((m^2+1)^{\frac12}+1))
       +\|\,|A|\,\|_{L^2(K)}\big]\,\\
  &&\qquad\qquad\qquad\qquad\qquad\qquad
  \times\big[\|\nabla \phi \|_{L^{\infty}(K)}+\|\phi \|_{L^{\infty}(K)}\big],
\end{eqnarray*}%
}
for all $\phi \in C_0^{\infty}({\Bbb R}^d)$ with $\text{supp}\,\phi \subseteq K$.
This says that $(H_{A,m})^{\alpha} u$ is a { continuous} linear functional
on $C_0^{\infty}({\Bbb R}^d)$, and so a distribution on ${\Bbb R}^d$.

\bigskip
Next, we study some properties of $(H_{A,m})^{\alpha}$ in the case
$A\equiv 0$, namely,
$(H_{0,m})^{\alpha} \equiv (-\Delta+m^2)^{\alpha/2},\,\, 0<\alpha \leq 1$.
This is the $\frac{\alpha}2$-power of the nonnegative selfadjoint operator
$H_{0,m} \equiv -\Delta+ m^2$
on $L^2({\Bbb R}^d)$ or also
a pseudo-differential operator defined through Fourier transform
having the symbol $(|\xi|^{2}+m^{2})^{\alpha/2}$. The function $\xi \mapsto
(|\xi|^{2}+m^{2})^{\alpha/2}-m^{\alpha}$ is conditionally negative definite
in ${\Bbb R}^d$
(e.g. [ReSi78, Appendix 2 to XIII.12, pp. 212--222]; [IkW81/89, p.65]),
so that, for each fixed  $t>0$, the function
$e^{-t[(|\xi|^{2}+m^{2})^{\alpha/2}-m^{\alpha}]}$ is positive definite.
 We note that this is a specific case of a Bernstein function,
providing the kinetic term of more general non-local Schr\"odinger operators
which we have studied in \cite{HILo12}.

As a result, its Fourier transform is a nonnegative function for each $t>0$,
which is nothing but the integral kernel $k_0^{m,\alpha}(t,x)$
of the semi-group $e^{-t[(H_{0,m})^{\alpha}-m^{\alpha}]}$ satisfying
$\int_{{\Bbb R}^d} k_0^{m,\alpha}(t,x) dx =1$.
We see further the operator $(H_{0,m})^{\alpha}u$, say with
$u \in C_0^{\infty}({\Bbb R}^d)$, have an integral operator representation:
\begin{align}
((H_{0,m})^{\alpha} u)(x)
&\equiv ([-\Delta+m^2]^{\frac{\alpha}2} u)(x)
 \equiv ({\cal F}^{-1}(|\xi|^{2}+m^{2})^{\frac{\alpha}2}{\cal F}u)(x)
                                                       \nonumber\\
&= m^{\alpha}u(x)- \int_{|y|>0} [u(x+y)-u(x)-I_{\{|y|<1\}}\,
      y\cdot \nabla_x u(x)]n^{m,\alpha}(dy),
\end{align}%
where $n^{m,\alpha}(dy)$ is a $\sigma$-finite measure on
${\Bbb R}^d \setminus \{0\}$
depending on $m\geq 0$ and $0< \alpha\leq 1$, called {\it L\'evy measure},
which satisfies $\int_{|y|>0}\frac{|y|^2}{1+|y|^2} n^{m,\alpha}(dy) < \infty$.
The L\'evy measure is known [IkW62, Example.1, p.81]
to be given from $k_0^{m,\alpha}(t,x)$ through
\begin{equation}
  \frac1{t} k_0^{m,\alpha}(t,dy) \quad\rightarrow\quad  n^{m,\alpha}(dy),
   \qquad t \downarrow 0.
\end{equation}%
In our case, it has density: $n^{m,\alpha}(dy) = n^{m,\alpha}(y)dy$.

For the expressions for the integral kernel $k_0^{m,\alpha}(t,x)$
of $e^{-t[(H_{0,m})^{\alpha}-m^{\alpha}]}$ and
the density (function) $n^{m,\alpha}(y)$, see Appendix A, (A.2).
For $\alpha=1$, they are explicitly given
(e.g. [I89, (2.4ab), (2.2ab), pp.268--269], [LLos01, 7.11 (11)]) as
\begin{eqnarray} 
k_0^{m,1}(t,x)
&=& \left\{
 \begin{array}{ll}
 2\big(\frac{m}{2\pi}\big)^{(d+1)/2}
 \frac{te^{mt}K_{(d+1)/2}(m(x^2+t^2)^{1/2})}{(x^2+t^2)^{(d+1)/4}},
                             &\quad m>0, \\
 \frac{\Gamma\big(\frac{d+1}{2}\big)}{\pi^{(d+1)/2}}
             \frac{t}{(x^2+t^2)^{(d+1)/2}}, &\quad   m=0;
\end{array}  \right. \\  
n^{m,1}(y)
&=& \left\{
  \begin{array}{ll}
   2\big(\frac{m}{2\pi}\big)^{(d+1)/2}\,
   \frac{K_{(d+1)/2}(m|y|)}{|y|^{(d+1)/2}}\,, \, &\qquad\qquad\quad m>0,\\
        \frac{\Gamma\big(\frac{d+1}{2}\big)}{\pi^{(d+1)/2}}
        \frac{1}{|y|^{d+1}}\,,      &\qquad\qquad\quad m=0,
 \end{array}\right.     
\end{eqnarray}
where $K_{\nu}(\tau)$ is the modified Bessel function of
the third kind of order $\nu$, which satisfies
$0< K_{\nu}(\tau)
\leq C\max\{\tau^{-\nu}, \tau^{-\frac12}\} e^{-\tau}, \, \tau >0$
with a constant $C>0$ when $\nu \geq \frac12$.

\medskip
For our later use, let us calculate the commutator
$\big[(H_{A,m})^2, \psi\big]$ 
 with $\psi \in C_0^{\infty}({\Bbb R}^d)$. Here for two operators
$U$ and $V$, their {\it commutator} is denoted by $[U, V] := UV - VU$.
We have
\begin{align}%
\big[(H_{A,m})^2, \psi\big]
&= (-i\nabla-A)^2\psi - \psi(-i\nabla-A)^2 \nonumber\\
&= (i\nabla + A)(i\nabla\psi) + (i\nabla \psi)(i\nabla + A) \nonumber\\
&= [(\Delta \psi) + 2(i\nabla + A)(i\nabla\psi)]
\,\,\text{\rm or}
 = [(-\Delta \psi) + 2(i\nabla\psi)(i\nabla + A)],
\end{align}
as quadratic forms,
 i.e. for suitable functions $u,\, v$  on ${\Bbb R}^d$,
{
\begin{align*}%
\big(u,\big[(H_{A,m})^2, \psi\big]v\big)
&= \big((i\nabla\psi)(i\nabla + A)u,v \big)
        + \big(u,(i\nabla \psi)(i\nabla + A)v\big) \nonumber\\
&= \big(u,(\Delta \psi)v\big) + 2\big(u,(i\nabla + A)(i\nabla\psi)v\big)\\
\,\,\text{\rm or}
&=\big(u,(-\Delta \psi)v\big) + 2\big(u,(i\nabla\psi)(i\nabla + A)v\big).
\end{align*}
}
Here note that
$[i\nabla +A,\psi]v =(i\nabla \psi)v$ as well as
$[i\nabla + A, (i\nabla \psi)]v = (-\Delta\psi)v$.
In fact, it holds more generally with two  ${\Bbb R}^d$-valued functions
 $A$ and $B$ that for a function $v$ in ${\Bbb R}^d$
\begin{align}%
&[(H_{A,m})^2\psi - \psi(H_{B,m})^2]v \nonumber\\
&= (i\nabla+A)\big((i\nabla \psi)+\psi A\big)v
    + \big((i\nabla\psi)-\psi B\big)(i\nabla +B)v
  + \psi A(i\nabla v)-i\nabla (\psi B v).
\end{align}

Indeed, the left-hand side of (2.11) can be seen to be equal to
\begin{align*}%
&\big[(-i\nabla-A)^2\psi - \psi(-i\nabla-B)^2]v\\
&= \big[(i\nabla+A)(i\nabla+A)\psi - \psi(i\nabla+B)(i\nabla+B)]v\\
&= (i\nabla+A)\big[\big((i\nabla \psi)+\psi A\big)v
        +\big(i\nabla(\psi v) -(i\nabla\psi)v\big)\big]\\
 &\qquad\qquad\qquad\qquad
   +\big[\big((i\nabla\psi)-\psi B\big)
       -\big(\psi(i\nabla)+ (i\nabla\psi)\big)\big](i\nabla +B)v\\
&= (i\nabla+A)\big((i\nabla \psi)+\psi A\big)v
       + (i\nabla+A)(\psi i\nabla v)\\
 &\qquad\qquad\qquad\qquad +\big((i\nabla\psi)-\psi B\big)(i\nabla +B)v
       -i\nabla\big(\psi(i\nabla +B)v\big) \\
&= (i\nabla+A)\big((i\nabla \psi)+\psi A\big)v
    + \big((i\nabla\psi)-\psi B\big)(i\nabla +B)v
    + \psi A(i\nabla v)-i\nabla (\psi B v).
\end{align*}%
This shows (2.11). Taking $B=A$ in (2.11) yields the third member of
(2.10), which implies the fourth and fifth members.

\bigskip
For the next lemma, we briefly mention the {\it weak} $L^1$-space
${L_w^1}(X)$, given a measurable subset $X$ of ${\Bbb R}^d$. It is
by definition the linear space of all measurable function $f$ on $X$
such that
\begin{equation}
 \|f\|_{L_w^1} := \sup_{a >0} a\, |\{x \in X;\, |f(x)| > a \}|
\end{equation}
is finite, where $|Y|$ denotes the volume (Lebesgue measure) of
the measurable set $Y \subseteq {\Bbb R}^d$.  ${L_w^1}(X)$ is not a
Banach space,
because $\|f\|_{L^1_w}$ is not a norm but a quasi-norm, as
it does not satisfy the triangle inequality. However, it holds that
$\|f+g\|_{L_w^1} \leq 2(\|f\|_{L_w^1} + \|g\|_{L_w^1})$.
It is shown that ${L_w^1}(X)$ is a {\it quasi-normed complete} linear space
(see e.g. [G10, Def.1.1.5, pp.5--6]).
We have $\|f\|_{L_w^1} \leq \|f\|_{L^1}$, so that
${L^1}(X) \subseteq {L_w^1}(X)$.
If $f_n \rightarrow f$ in $L_w^1$, then the $\{f_n\}$ converges to $f$
in measure (e.g. [G10, Prop.1.1.9, p.7]).
We say ``{\it $f$ is locally in $L_w^1$}", if for every compact set $K$
in ${\Bbb R}^d$, $f$ belongs to $L_w^1(K)$.
In some literatures ${L_w^1}(X)$ may be denoted also by
$L^{1,\infty}(X)$ ({\it Lorentz space}).

\begin{lemma}
Let $0< \alpha \leq 1$. Let $\psi \in C_0^{\infty}({\Bbb R}^d)$.
Then for the commutator $[(H_{0,m})^{\alpha},\psi]$, it holds,
 with a constant $C_{\alpha}$ dependent
on $\psi$ and $\alpha$ but independent of $m\geq 0$,
that
(i) for $1< p<\infty$,
\begin{equation}%
 \|[(H_{0,m})^{\alpha}, \psi]u\|_{L^p}
= \|(H_{0,m})^{\alpha}(\psi u) -\psi (H_{0,m})^{\alpha} u\|_{L^p}
  \leq C_{\alpha}\|u\|_{L^p},
\end{equation}
for all $u \in L^p({\Bbb R}^d)$. Therefore if both
$u$ and $(H_{0,m})^{\alpha} (\psi u)$ are in $L^p$, then
$\psi (H_{0,m})^{\alpha} u$ is in $L^p$, and
$$
 \|\psi (H_{0,m})^{\alpha} u\|_{L^p} \leq C_{\alpha}\|u\|_{L^p}
                 + \|(H_{0,m})^{\alpha} (\psi u)\|_{L^p};
$$
(ii) for $p=1$,
\begin{equation}
 \|[(H_{0,m})^{\alpha}, \psi]u\|_{L^1_w}
= \|(H_{0,m})^{\alpha}(\psi u) -\psi (H_{0,m})^{\alpha} u\|_{L^1_w}
  \leq C_{\alpha}\|u\|_{L^1},
\end{equation}%
for all $u \in L^1({\Bbb R}^d)$.
\end{lemma}

\bigskip\noindent
{\it Remark}. Inequality (2.13) does not hold for $p=1$, and
instead we have (2.14) with the $L^1$-norm on the left-hand side
replaced by the $L_w^1$-quasi-norm. This is dependent on the
Calder\'on--Zygmund theorem (For this see Proposition 2.4 below).

\bigskip
{\it Proof of Lemma 2.3}.
(i) As the second-half assertion follows from the first, i.e.
inequality (2.13), we have only to show (2.13), and even
only for $u \in C_0^{\infty}({\Bbb R}^d)$, since $C_0^{\infty}({\Bbb R}^d)$
is dense in $L^2({\Bbb R}^d)$.
The proof for the case $\alpha=1$ was given in [ITs92, p.274, Lemma 2.3] by
using the integral operator representation (2.6) of
$H_{0,m} = \sqrt{-\Delta+m^2}$.
The proof for the case $0<\alpha < 1$ is similar. So we only give an
outline.

Use (2.6) to rewrite $[(H_{0,m})^{\alpha},\psi]$ as
\begin{eqnarray}%
([(H_{0,m})^{\alpha},\psi]u)(x)
&=\!& -\int_{|y|>0} [\psi(x+y)-\psi(x)-I_{\{|y|<1\}}
      y\cdot \nabla_x \psi(x)]u(x+y)n^{m,\alpha}(dy)\nonumber\\
 &&\qquad
   -\int_{0<|y|<1} y\cdot \nabla_x\psi(x)[u(x+y)-u(x)]n^{m,\alpha}(dy)
                                                    \nonumber\\
&=\!&: (I_1u)(x) + (I_2u)(x).
\end{eqnarray}
We estimate the $L^p$ norms of $I_1u$ and $I_2u$ in the last member.

First, rewrite $I_1u$ as
\begin{eqnarray*}%
(I_1u)(x)
&=& -\int_{0<|y|<1} [\psi(x+y)-\psi(x)-I_{\{|y|<1\}}
      y\cdot \nabla_x \psi(x)]u(x+y)n^{m,\alpha}(dy)\nonumber\\
 &&\qquad -\int_{|y|\geq 1} [\psi(x+y)-\psi(x)]u(x+y)n^{m,\alpha}(dy).
\end{eqnarray*}%
Hence
$$
|(I_1u)(x)|\leq \|\nabla^2\psi\|_{L^{\infty}}
             \int_{0<|y|<1} |y|^2 |u(x+y)|n^{m,\alpha}(dy)
     +2\|\psi\|_{L^{\infty}} \int_{|y|>1} |u(x+y)|n^{m,\alpha}(dy),
$$
so that for $1\leq p< \infty$
$$
\|I_1 u\|_{L^p}
= \Big(\int |(I_1u)(x)|^p dx\Big)^{\frac1{p}}
\leq\big(n_1^{m,\alpha} \|\nabla^2\psi\|_{L^{\infty}}
     + 2n_{\infty}^{m,\alpha}\|\psi\|_{L^{\infty}}\big) \|u\|_{L^p},
$$
where
\begin{equation}
n_{\infty}^{m,\alpha} := \int_{|y|\geq 1} n^{m,\alpha}(dy), \qquad
n_{\kappa}^{m,\alpha} := \int_{0<|y|<1} |y|^{1+\kappa}n^{m,\alpha}(dy)\,,
\end{equation}%
where the former is finite, and the latter is finite for all $0<\kappa \leq 1$.

Next, for $I_2u$ we use the following known fact for an operator $T$
on $L^p({\Bbb R}^d)$ with Calder\'on--Zygmund kernel
$K :{\Bbb R}^d\setminus \{0\} \rightarrow {\Bbb C}$
(e.g. [St70, II.3, pp.35--42], [G10, Theorem 5.3.3, p.359],
[MSc13, Def.7.1, Prop.7.4, Theorem 7.5, pp.166--172]).
It is the integral kernel which satisfies, for some constant $B>0$,
the following conditions:

\medskip
(i) $|K(x)|\leq B|x|^{-d}$ for all $x \in {\Bbb R}^d$; $\,\,$

(ii) $\int_{|x|\geq 2|y|} |K(x)-K(x-y)|\,dx \leq B$ for all $y \not= 0$; $\,\,$

(iii) $\int_{R_1<|x|<R_2} K(x)\ dx =0$ for all $0<R_1<R_2<\infty$.

\bigskip\noindent
\begin{prp}
Let
$$
(Tf)(x) := \lim_{\varepsilon\downarrow 0}
  \int_{|x-y|\geq \epsilon} K(x-y)f(y)dy.
$$
Then
\begin{align*}
 \|Tf\|_{L^p} &\leq C_p\|f\|_{L^p}\,,  &1<p<\infty,\\
 \|Tf\|_{L_w^1}
  &\equiv \sup_{a>0}\, a\,\big|\{x\in {\Bbb R}^d;\, |(Tf)(x)|>a\}\big|
  \leq C_1\|f\|_{L^1}\,,               &p=1.
\end{align*}
\end{prp}

This proposition is going to be used in the proof of Lemma 2.3 (i).

\bigskip
We continue the proof of Lemma 2.3 (i). It still remains to deal with
$I_2u$, which is rewritten as
$$
 (I_2u)(x) = -\sum_{j=1}^d \lim_{\varepsilon\downarrow 0}
           \int_{\varepsilon\leq |y|<1}\,\partial_{x_j}\psi(x)\,
           (x_j-y_j) n^{m,\alpha}(x-y)u(y)dy.
$$
Here each $x_j\cdot n^{m,\alpha}(x),\, 1\leq j \leq d$, is a
Calder\'on--Zygmund kernel (see Appendix A, (A.2)), so that
we have by Proposition 2.4 with $1 \leq p<\infty$
there exists a constant $C_p>0$ such that
\begin{eqnarray*}%
\|I_2u\|_{L^p}
&\leq& C_p \|\nabla \psi\|_{L^{\infty}} \|u\|_{L^p}, \quad 1<p<\infty,\\
\|I_2u\|_{L^1_w}
&=& \text{\rm sup}_{a>0}\,\,
   a\big|\{x \in {\Bbb R}^d;\; |(I_2u)(x)| > a\}\big|
\leq C_1 \|\nabla \psi\|_{L^{\infty}} \|u\|_{L^1},  \quad p=1.
\end{eqnarray*}%
Thus we obtain
\begin{eqnarray*}%
\|[(H_{0,m})^{\alpha},\psi]u\|_{L^p}
&\leq& \|I_1 u\|_{L^p} + \|I_2 u\|_{L^p}\\
&\leq& \big(n_1^{m,\alpha}\|\nabla^2 \psi\|_{L^{\infty}}
    +2n_{\infty}^{m,\alpha}\|\psi\|_{L^{\infty}}
    +C_p\|\nabla \psi\|_{L^{\infty}}\big) \big\|u\|_{L^p}, \\
\end{eqnarray*}%
showing (i) for $1<p<\infty$.

Next, for (ii) for $p=1$, we have
\begin{eqnarray*}%
\|[(H_{0,m})^{\alpha},\psi]u\|_{L^1_w}
&\leq& 2\big(\|I_1 u\|_{L^1_w} + \|I_2 u\|_{L^1_w}\big)
 \leq 2\|I_1 u\|_{L^1} + 2\|I_2 u\|_{L^1_w}\\
&\leq& 2\big(n_1^{m,\alpha}\|\nabla^2 \psi\|_{L^{\infty}}
    +2n_{\infty}^{m,\alpha}\|\psi\|_{L^{\infty}}
      +C_1\|\nabla \psi\|_{L^{\infty}}\big) \big\|u\|_{L^1},
\end{eqnarray*}%
because $\|I_1u\|_{L^1_w}\leq \|I_1u\|_{L^1}$.
This shows (ii), ending the proof of Lemma 2.3.
\qed

\bigskip
When $A \in L^2_{\text{\rm \rm loc}}$, our selfadjoint operator
$S:= (-i\nabla -A)^2+m^2$ originally is being defined
as the selfadjoint operator in $L^2({\Bbb R}^d)$
associated with the closed quadratic form (2.1).
As already noted in the proof of Lemma 2.1 (i), it also makes sense as
an operator in the spaces $L^p({\Bbb R}^d), \,\ 1 \leq p < \infty$,
referring to the result [Si79, Theorem 2.3] or [Si82, Sect. B13]) that
the Schr\"odinger semigroup $e^{-tS}= e^{-t[(-i\nabla -A)^2+m^2]}$
satisfies
\begin{equation}
 |e^{-t[(-i\nabla -A)^2+m^2]}g| \leq e^{-t[-\Delta+m^2]}|g|
\end{equation}
pointwise for any $g \in L^2({\Bbb R}^d)$.
This yields that for $1 \leq p <\infty$, $e^{-t(H_{A,m})^2}$ is a bounded
operator of $L^p({\Bbb R}^d)$ into itself for all $t>0$,
which also is a contraction semigroup.

Thus, the fractional powers of $S$ such as
$S^{\frac{\alpha}2} = (H_{A,m})^{\alpha}$
in (2.3) equally make sense in $L^p({\Bbb R}^d)$.

\bigskip
Now, we give two crucial lemmas, Lemmas 2.5 and 2.6.

\begin{lemma}
Let $0<\alpha < 1$ and assume that $A \in [L^2_{\text{\rm loc}}({\Bbb R}^d)]^d$.
 Then: (i) if $u \in L^2({\Bbb R}^d)$, one has
for $\chi,\,\, \psi \in C_0^{\infty}({\Bbb R}^d)$
\begin{align}%
\|\chi [(H_{0,m})^{\alpha}\psi -\psi (H_{A,m})^{\alpha}]u\|_{L^1}
&\equiv \big\|\chi \big([-\Delta+m^2]^{\frac{\alpha}2}\psi
 - \psi[(-i\nabla-A)^2+m^2]^{\frac{\alpha}2}\big)u\big\|_{L^1}\nonumber\\
&\leq C_{\alpha,A,m,\chi, \psi} \|u\|_{L^2},    
\end{align}%
where $C_{\alpha,A,\chi, \psi}$ is a constant which depends
on $0<\alpha<1$, $A$, $m>0$, $\chi$ and $\psi$, and which tends
to $\infty$ as $\alpha \uparrow 1$.

(ii) In particular, when $A=0$, (2.18) reads:
if $u \in L^2({\Bbb R}^d)$, one has
\begin{equation}
\|\chi [(H_{0,m})^{\alpha}, \psi]u\|_{L^1}
\leq C_{\alpha,0,m,\chi,\psi} \|u\|_{L^2}.
\end{equation}
\end{lemma}

\bigskip
For $A=0$, inequality (2.19) appears more useful in comparison with (2.14),
Lemma 2.3.

\bigskip
{\it Proof of Lemma 2.5}.

\medskip
(i) We have only to show (2.18) when $u \in C_0^{\infty}({\Bbb R}^d)$,
since $C_0^{\infty}({\Bbb R}^d)$ is dense in $L^2({\Bbb R}^d)$.
Note then that
$H_{0,m}u$ and $H_{A,m}u $ belong to $L^2({\Bbb R}^d)$.

We use formula (2.3) for $(H_{0,m})^{\alpha}$ as well as $(H_{A,m})^{\alpha}$
to calculate
\begin{eqnarray*}%
&&[(H_{0,m})^{\alpha}\psi - \psi(H_{A,m})^{\alpha}]u \\
&=& \frac1{\Gamma(\frac{2-\alpha}2)}\int_0^{\infty} t^{-\frac{\alpha}2}
 \Big[e^{-t(H_{0,m})^2}(H_{0,m})^2\psi - \psi(H_{A,m})^2e^{-t(H_{A,m})^2}\Big] u\, dt\\
&=& \frac1{\Gamma(\frac{2-\alpha}2)}\int_0^{\infty}dt\, t^{-\frac{\alpha}2}
  \big(-\frac{d}{dt}\big)\Biggl[e^{-\theta t(H_{0,m})^2}\psi e^{-(1-\theta)t(H_{A,m})^2}
                                 \Biggr]_{\theta=0}^{\theta=1}u\\
&=&-\frac1{\Gamma(\frac{2-\alpha}2)}\int_0^{\infty}dt\, t^{-\frac{\alpha}2}
  \frac{d}{dt}\int_0^1 d\theta\, \frac{d}{d\theta}
  \Big[e^{-\theta t(H_{0,m})^2}\psi e^{-(1-\theta)t(H_{A,m})^2}\Big]u\\
&=& \frac1{\Gamma(\frac{2-\alpha}2)}\int_0^{\infty}dt\, t^{-\frac{\alpha}2}\\
 &&\quad \times
  \frac{d}{dt} \,\Big(t\int_0^1 d\theta\,
  \Big[e^{-\theta t(H_{0,m})^2}\big[(H_{0,m})^2\psi - \psi(H_{A,m})^2\big]
    e^{-(1-\theta)t(H_{A,m})^2}\Big]u\Big).
\end{eqnarray*}%
Then by integration by parts,
\begin{align}%
&[(H_{0,m})^{\alpha}\psi - \psi(H_{A,m})^{\alpha}]u\nonumber\\
&= \frac1{\Gamma(\frac{2-\alpha}2)}\Biggl[t^{-\frac{\alpha}2+1} \int_0^1
  \Big(e^{-\theta t(H_{0,m})^2}\big[(H_{0,m})^2\psi - \psi(H_{A,m})^2\big]
    e^{-(1-\theta)t(H_{A,m})^2}\Big)u\, d\theta\Biggr]_{t=0}^{t=\infty}
                                                              \nonumber\\
 &\,\,+\frac{\alpha}{2\Gamma(\frac{2-\alpha}2)}\int_0^{\infty}dt\,
  t^{-\frac{\alpha}2} \int_0^1 d\theta\,
  \Big(e^{-\theta t(H_{0,m})^2}\big[(H_{0,m})^2\psi - \psi(H_{A,m})^2\big]
    e^{-(1-\theta)t(H_{A,m})^2}\Big)u                         \nonumber\\
&= \frac{\alpha}{2\Gamma(\frac{2-\alpha}2)}\int_0^{\infty}dt\,
  t^{-\frac{\alpha}2} \int_0^1 d\theta\,
  \Big(e^{-\theta t(H_{0,m})^2}\big[(H_{0,m})^2\psi - \psi(H_{A,m})^2\big]
    e^{-(1-\theta)t(H_{A,m})^2}\Big)u.  
\end{align}
Here we make two observations related to (2.20).
First for its second member, the boundary value at
$t \rightarrow \infty$
of the first term also vanishes, because the part
$$
e^{-\theta t(H_{0,m})^2}\big[\cdots \big]e^{-(1-\theta)t(H_{A,m})^2}
=e^{-\theta t(-\Delta+m^2)} \big[\cdots \big]
                e^{-(1-\theta)t[(-i\nabla-A)^2+m^2]}
$$
contains the factor  $e^{-m^2t}$.
Second for its last member, note that the middle factor in the integrand
is, by (2.11) with $A :=0,\,\, B := A$, equal to
\begin{equation}
[(H_{0,m})^2\psi - \psi(H_{A,m})^2]
= \big[i\nabla\big((i\nabla \psi)-\psi A\big)
   +\big((i\nabla\psi)-\psi A\big)(i\nabla +A)\big]
\end{equation}
as quadratic forms.

Substituting (2.21) into (2.20), we have with
$\chi \in C_0^{\infty}({\Bbb R}^d)$
\begin{align}%
&\chi [(H_{0,m})^{\alpha}\psi - \psi(H_{A,m})^{\alpha}]u \nonumber\\
&= \frac{\alpha}{2\Gamma(\frac{2-\alpha}2)} \int_0^{\infty}dt\,
  t^{-\frac{\alpha}2}\int_0^1 d\theta\,
 \chi \Big(e^{-\theta t(H_{0,m})^2}i\nabla\big((i\nabla \psi)-\psi A\big)
  \, e^{-(1-\theta)t(H_{A,m})^2}\Big)u \nonumber\\
&\quad +\frac{\alpha}{2\Gamma(\frac{2-\alpha}2)}\int_0^{\infty}dt\,
  t^{-\frac{\alpha}2}\int_0^1 d\theta\,
 \chi \Big(e^{-\theta t(H_{0,m})^2} \big((i\nabla\psi)- \psi A\big)(i\nabla +A)
  \, e^{-(1-\theta)t(H_{A,m})^2}\Big)u  \nonumber\\
&=: I_3u + I_4u\,.    
\end{align}%
We estimate the $L^1$ norm for $I_3u$ and $I_4u$ in (2.22).
 Note that $e^{-t(-i\nabla-A)^2},\,\, t\geq 0$, is a contraction on
$L^p({\Bbb R}^d)$, $1\leq p \leq \infty$.

\medskip
First, for $I_3u$, integrate its absolute value in $x$ to get
\begin{eqnarray}%
\|I_3u\|_{L^1}
&\leq& \frac{\alpha}{2\Gamma(\frac{2-\alpha}2)}\int_0^{\infty}
        t^{-\frac{\alpha}2}e^{-m^2t}dt
        \int_0^1 d\theta \nonumber\\
 &&\,\times
   \big\|\chi [e^{-\theta t(-\Delta)} (i\nabla)]
   \big((i\nabla \psi)-\psi A\big)\,
   e^{-(1-\theta)t(-i\nabla-A)^2}u\big\|_{L^1}\,.
\end{eqnarray}
Then by Lemma 2.1 (ii) for $p=1$, the Schwarz inequality
and Lemma 2.1 (i)
{\begin{align*}%
\|I_3u\|_{L^1}
&\leq \frac{\alpha}{2\Gamma(\frac{2-\alpha}2)}\int_0^{\infty}
        t^{-\frac{\alpha}2}e^{-m^2t}dt
        \int_0^1 d\theta \nonumber\\
 &\,\times
   \|\chi\|_{L^{\infty}}
   \|e^{-\theta t(-\Delta)} (i\nabla)\|_{[L^1]^d\to [L^1]^d}\,
   \|\big((i\nabla \psi)-\psi A\big)\,
    e^{-(1-\theta)t(-i\nabla-A)^2}u\big\|_{L^1\to [L^1]^d}\nonumber\\
&\leq \frac{\alpha}{2\Gamma(\frac{2-\alpha}2)}\int_0^{\infty}
        t^{-\frac{\alpha}2}e^{-m^2t}dt
        \int_0^1 d\theta \nonumber\\
 &\,\times
   \|\chi\|_{L^{\infty}} C_{11} (\theta t)^{-1/2}
   \|(i\nabla \psi)-\psi A\|_{L^2}\,
   \|e^{-(1-\theta)t(-i\nabla-A)^2}u\|_{L^2}\nonumber\\
&\leq \frac{C_{11}\alpha}{2\Gamma(\frac{2-\alpha}2)}\int_0^{\infty}
        t^{-\frac{1+\alpha}2}e^{-m^2t}dt
        \int_0^1 \tfrac{d\theta}{\theta^{1/2}}
  \|\chi\|_{L^{\infty}} \|(i\nabla \psi)-\psi A\|_{L^2}\|u\|_{L^2}
                                                     \nonumber\\
\end{align*}%
}
Here recall
that $\|(i\nabla \psi)-\psi A\big\|_{L^2} <\infty$ by assumption on $A$ and
notice also that
$$\int_0^{\infty}t^{-\frac{1+\alpha}2} e^{-m^2t}dt
  = \Gamma(\tfrac{1-\alpha}2)\,m^{-\frac{1-\alpha}2},
$$
which diverges as $\alpha \uparrow 1$ with $m>0$.
Thus we have
\begin{eqnarray}%
\|I_3u\|_{L^1}
&\leq& \frac{C_{11}\alpha\Gamma(\frac{1-\alpha}2)}
        {\Gamma(\frac{2-\alpha}2)\,m^{\frac{1-\alpha}2}}\,
   \big\|(i\nabla \psi)-\psi A\big\|_{L^2}
   \|\chi\|_{L^{\infty}}\big\|u\big\|_{L^2}.
\end{eqnarray}

Next for $I_4u$, in a similar way, we have from (2.22)
\begin{eqnarray}%
\|I_4u\|_{L^1}
&\leq& \frac{\alpha}{2\Gamma(\frac{2-\alpha}2)} \int_0^{\infty}
  t^{-\frac{\alpha}2} e^{-m^2t}dt
  \int_0^1 d\theta
  \big\|\chi \big(e^{-\theta t(-\Delta)}\big((i\nabla\psi)- \psi A\big)
                                                          \nonumber\\
   &&\qquad\qquad\qquad\qquad \times
   [(i\nabla+A)\, e^{-(1-\theta)t(-i\nabla-A)^2}]\big)u\|_{L^1}
\end{eqnarray}
Then by the Schwarz inequality and Lemma 2.1 (iv)
\begin{align*}%
\|I_4u\|_{L^1}
&\leq \frac{\alpha}{2\Gamma(\frac{2-\alpha}2)} \int_0^{\infty}
  t^{-\frac{\alpha}2} e^{-m^2t}dt
  \int_0^1 d\theta
  \|\chi\|_{L^{\infty}} \|e^{-\theta t(-\Delta)}\|_{L^1\to L^1}
                                                          \nonumber\\
 &\qquad\qquad\qquad\qquad \times
  \|\big((i\nabla\psi)- \psi A\big)
  [(i\nabla+A)\, e^{-(1-\theta)t(-i\nabla-A)^2}]\big)u\|_{L^1}\\
&\leq \frac{\alpha}{2\Gamma(\frac{2-\alpha}2)} \int_0^{\infty}
  t^{-\frac{\alpha}2} e^{-m^2t}dt
  \int_0^1 d\theta
  \|\chi\|_{L^{\infty}} \|e^{-\theta t(-\Delta)}\|_{L^1\to L^1}
  \|(i\nabla\psi)- \psi A\|_{L^2}
                                                          \nonumber\\
   &\qquad\qquad\qquad\qquad \times
   \|(i\nabla+A)\, e^{-(1-\theta)t(-i\nabla-A)^2}]u\|_{L^2}\\
&\leq \frac{\alpha}{2\Gamma(\frac{2-\alpha}2)} \int_0^{\infty}
  t^{-\frac{\alpha}2} e^{-m^2t}dt
  \int_0^1 d\theta \|\chi\|_{L^{\infty}}  \|(i\nabla\psi)- \psi A\|_{L^2}
   \big(\tfrac{d}{2e(1-\theta)}\big)^{1/2}\|u\|_{L^2}\\
&=  \big(\frac{d}{2e}\big)^{1/2}\frac{\alpha}{2\Gamma(\frac{2-\alpha}2)}
  \int_0^{\infty}  t^{-\frac{1+\alpha}2} e^{-m^2t}dt
  \int_0^1 \tfrac{d\theta}{(1-\theta)^{1/2}} \|\chi\|_{L^{\infty}}
       \|(i\nabla\psi)- \psi A\|_{L^2}\|u\|_{L^2}.
\end{align*}%
Then we have
\begin{equation}%
\|I_4u\|_{L^1}
\leq \big(\frac{d}{2e}\big)^{1/2}\frac{\alpha\Gamma(\frac{1-\alpha}2)}
  {\Gamma(\frac{2-\alpha}2)\,m^{\frac{1-\alpha}2}} \,
  \big\|(i\nabla \psi) -\psi A\big\|_{L^2} \|\chi\|_{L^{\infty}}
  \|\chi\|_{L^{\infty}}\|u\|_{L^2}.
\end{equation}

\bigskip

Putting (2.24) and (2.26) together in view of (2.22),
we have
\begin{align}%
&\|\chi [(H_{0,m})^{\alpha}\psi -\psi(H_{A,m})^{\alpha}]u\|_{L^1}\nonumber\\
&\leq 2\big(\|I_3u\|_{L^1} + \|I_4u\|_{L^1}\big) \nonumber\\
&\leq 2\frac{C_{11}+\big(\frac{d}{2e}\big)^{1/2}}{m^{\frac{1-\alpha}2}}
      \frac{\alpha\Gamma(\frac{1-\alpha}2)}
      {\Gamma(\frac{2-\alpha}2)} \,
      \big\|(i\nabla \psi)-\psi A\big\|_{L^2}
      \|\chi\|_{L^{\infty}} \big\|u\big\|_{L^2}.
\end{align}%
This yields (2.18), showing Lemma 2.5 (i).

(ii) Inequality (2.19) is immediately derived by putting $A=0$ in (2.18).

This shows Lemma 2.5 (ii), completing the proof of Lemma 2.5.
\qed

\bigskip
From Lemma 2.5 we have the following result which we shall need,
in particular, assertion (ii), in the proof of Theorem 1.1.

\begin{lemma}
Let $0<\alpha < 1$.  Assume
that $A \in [L^{2}_{\operatorname{\rm loc}}({\Bbb R}^d)]^d$.

(i) If $u  \in (C^{\infty}\cap L^2)({\Bbb R}^d)$, then
$(H_{A,m})^{\alpha} u $ is locally in $L^1({\Bbb R}^d)$.

(ii) If $u \in L^2({\Bbb R}^d)$ with
  $(H_{A,m})^{\alpha} u \in L^1_{\text{\rm loc}}({\Bbb R}^d)$,
then $(H_{0,m})^{\alpha} u$ is locally in $L^1({\Bbb R}^d)$.
\end{lemma}

\bigskip
{\it Proof}.
(i) Let $u \in (C^{\infty}\cap L^2)({\Bbb R}^d)$. Then
for $\psi \in C_0^{\infty}({\Bbb R}^d)$,
$$
\psi (H_{A,m})^{\alpha} u = (H_{0,m})^{\alpha}(\psi u)
      +  \big(\psi (H_{A,m})^{\alpha} - (H_{0,m})^{\alpha} \psi \big)u.
$$
Put $K = \text{supp}\,\psi$.
Then, since $\psi u$ is in $C_0^{\infty}({\Bbb R}^d)$,
the first term $(H_{0,m})^{\alpha}(\psi u)$ on the right-hand side
belongs to $L^2({\Bbb R}^d)$, as we can see
from (2.6) (with $\psi u$ instead of $u$) or
Lemma 2.2 (2.4) with $A=0$ (with $\psi u$ instead of $\varphi$).
For the second term restricted to $K$, it belongs to $L^1(K)$,
as we see by Lemma 2.5 (2.18).
Therefore $\psi (H_{A,m})^{\alpha} u$ is in $L^1(K)$, so that
$(H_{A,m})^{\alpha} u$ is locally in $L^1({\Bbb R}^d)$.
This proves the assertion (i).

\medskip
(ii) Let $u \in L^2$ with $(H_{A,m})^{\alpha} u \in L^1_{\text{\rm loc}}$ and
let $K$ be an arbitrary compact subset of ${\Bbb R}^d$. Take
$\chi,\, \psi \in C_0^{\infty}({\Bbb R}^d)$ with $0 \leq \chi(x) \leq 1$
such that $\chi (x) = \psi(x) = 1$ on $K$.
Then since
$$%
\psi (H_{0,m})^{\alpha}u -\psi (H_{A,m})^{\alpha} u
= -[(H_{0,m})^{\alpha}, \psi] u
          + \big((H_{0,m})^{\alpha} \psi - \psi (H_{A,m})^{\alpha} \big) u,
$$
we have by Lemma 2.5 (2.18) with $A=0$ as well as with non-zero $A$
\begin{align*}%
&\|(H_{0,m})^{\alpha} u -\psi (H_{A,m})^{\alpha} u \|_{L^1(K)} \\
&= \|\chi\psi [(H_{0,m})^{\alpha} u - (H_{A,m})^{\alpha} u] \|_{L^1(K)}\\
&\leq \|\chi[(H_{0,m})^{\alpha},\psi] u\|_{L^1}
 +\|\chi \big((H_{0,m})^{\alpha}\psi - \psi(H_{A,m})^{\alpha} \big) u\|_{L^1}\\
&\leq (C_{\alpha,0,m,\chi,\psi} + C_{\alpha,A,m,\chi,\psi}) \|u\|_{L^2} < \infty\,.
\end{align*}%

Since by assumption $(H_{A,m})^{\alpha} u$ is locally in
$L^1({\Bbb R}^d)$, we have $(H_{0,m})^{\alpha} u$ is locally $L^1({\Bbb R}^d)$.
This proves the assertion (ii), ending the proof of Lemma 2.6.
\qed

\section{Proof of Theorems} 
We show only Theorem 1.1 and Theorem 1.2.
As for Theorem 1.3,
 essential selfadjointness of $H_{A,V,m}$ follows from Theorem 1.1 by
its standard application in Kato's original paper [K72].
In fact, one can show in the same way as in [I89, Theorem 5.1].
So the proof is omitted. The assertion
that $H_{A,V,m} = H_{A,m}+V \geq m$ is trivial because $H_{A,m} \geq m$.

\medskip
In this section, we keep assuming that $m>0$ before come to the final part
(iii) of the proof of Theorem 1.1.

\subsection{Proof of Theorem 1.1}

The proof will proceed similarly to Kato's original proof [K72]
 (e.g. [{ReSi75}, Theorems X.27 (p.183), X.33 (p.188)]) for the
magnetic {\it nonrelativistic} Schr\"odinger operator
$\frac1{2m}(-i\nabla -A(x))^2$ and to a modified one [I89], [ITs92]
for {\it another} magnetic relativistic Schr\"odinger operator.
However, if one could show the assumption of the theorem
that $u \in L^2$ with $H_{A,m}u \in L^1_{\text{\rm loc}}$ implies that
{ $\partial_j u \in L^1_{\text{\rm loc}}\,, \, 1\leq j \leq d$, and/or
$H_{0,m}u \in L^1_{\text{\rm loc}}\,$,
}
there should be no problem.
The obstruction seems to come from the fact that the operators
$\partial_j\cdot (-\Delta+m^2)^{-1/2}$, $1\leq j \leq d$,
 are {\it not} bounded from $L^1$ to $L^1$,
though they are {\it bounded} from $L^1$ to $L^1_{w}$.
The strategy we adopt to cope with this difficulty is, in the beginning,
to make a detour by considering the case $(H_{A,m})^{\alpha}$ for
$\alpha<1$, leaving the very case $\alpha=1$ aside, however,
to handle the local convergence in $L^1$.
In fact, in the first stage  (Lemmas 3.1 and 3.2),
we show first that if $(H_{A,m})^{\alpha}u \in L^1_{\text{\rm loc}}$, then
$(H_{A,m})^{\alpha}u^{\delta} \rightarrow (H_{A,m})^{\alpha}u$
locally in $L^1$ as $\delta \downarrow 0$,
and making use of  Lemma 2.6 saying that
$(H_{0,m})^{\alpha}u$ is locally in $L^1$. 
Next we show that
the assumption $H_{A,m}u \in L^1_{\text{\rm loc}}$ implies that
$(H_{A,m})^{\alpha}u \in L^1_{\text{\rm loc}}$ for $0<\alpha<1$,
and $(H_{A,m})^{\alpha}u$ converges to $H_{A,m}u$ in
$L^1_{\text{\rm loc}}$ as $\alpha \uparrow 1$.
In the second, main stage, with $m>0$,  we show {\it first for}
 $0<\alpha<1$ that the asserted inequality, i.e.
\begin{equation}%
  \text{Re} ((\text{sgn}\, u)[(H_{A,m})^{\alpha} -m^{\alpha}]u)
                    \geq [(H_{0,m})^{\alpha}-m^{\alpha}]|u|,
\end{equation}
holds, and {\it next for} $\alpha=1$, using the just above mentioned
fact that $(H_{A,m})^{\alpha}u \rightarrow H_{A,m}u$ in
$L^1_{\text{\rm loc}}$ as $\alpha \uparrow 1$.
The final stage will deal with the remaining case for $m=0$ and $\alpha=1$.

\bigskip
We provide two lemmas playing a crucial role in the proof of Theorem 1.1.

\medskip
For a function $f$ locally in $L^1({\Bbb R}^d)$,
we write its {\it mollifier} as
$f^{\delta} = \rho_{\delta}*f, \,\, 0<\delta \leq 1$, where
$\rho_{\delta}(x) := \delta^{-d}\rho(x/\delta)$, and $\rho(x)$ is a nonnegative
$C^{\infty}$ function ${\Bbb R}^d$ with compact support
$\text{supp}\, \rho \subseteq \{x;\, |x|\leq 1\}$ and
$\int \rho(x) dx =1$.

\begin{lemma} 
Let $0<\alpha< 1$. Let $u \in L^2({\Bbb R}^d)$, so that
$u^{\delta} := \rho_{\delta}* u \rightarrow u$ in $L^2$ as
$\delta \downarrow 0$.
If $(H_{A,m})^{\alpha} u \in L^1_{\text{\rm loc}}({\Bbb R}^d)$,
then $(H_{A,m})^{\alpha} u^{\delta} = [(-i\nabla -A)^2+m^2]^{\alpha/2}u^{\delta}
  \rightarrow  (H_{A,m})^{\alpha} u =[(-i\nabla -A)^2+m^2]^{\alpha/2}u$
locally in $L^1({\Bbb R}^d)$ as $\delta \downarrow 0$.
\end{lemma}

\bigskip
{\it Proof}.
 Let $u \in L^2$ and 
$(H_{A,m})^{\alpha}u \in L^{1}_{\text{\rm loc}}({\Bbb R}^d)$.
Then by Lemma 2.6 (ii), $(H_{0,m})^{\alpha} u$ is locally in $L^1$
and since $u^{\delta} \in C^{\infty}\cap L^2$, we have by Lemma 2.6 (i)
that $(H_{A,m})^{\alpha} u^{\delta}$ is locally in $L^1$.
The important is:
thanks to the integral operator representation (2.6) of the operator
$(H_{0,m})^{\alpha}$, the convolution commutes with $(H_{0,m})^{\alpha}$.
Therefore we have
$((H_{0,m})^{\alpha} u)^{\delta} = (H_{0,m})^{\alpha} u^{\delta}$,
which converges to $(H_{0,m})^{\alpha} u$ locally in $L^1$
as $\delta \downarrow 0$.
Then for $K$ a compact subset in ${\Bbb R}^d$, let
$\chi, \, \psi \in C_0^{\infty}({\Bbb R}^d)$ with
$0\leq \chi(x)\leq 1$ on ${\Bbb R}^d$ and $\chi(x) = \psi(x) =1$ on $K$.
We have
\begin{align*}%
&\|(H_{A,m})^{\alpha} u^{\delta} - (H_{A,m})^{\alpha} u\|_{L^1(K)}\\
&= \|\chi \psi (H_{A,m})^{\alpha}(u^{\delta} - u)\|_{L^1(K)}\\
&= \|\chi \big[-(H_{0,m})^{\alpha}\psi +
  \big((H_{0,m})^{\alpha}\psi - \psi(H_{A,m})^{\alpha}\big)\big]
                                (u^{\delta} - u)\|_{L^1(K)}\\
&\leq \|\chi (H_{0,m})^{\alpha}\psi (u^{\delta} - u)\|_{L^1}
      +\|\chi [(H_{0,m})^{\alpha}\psi -\psi (H_{A,m})^{\alpha}]
                (u^{\delta} - u)\|_{L^1}.
\end{align*}%
The second term in the last member of the above inequality is,
by Lemma 2.5 (2.18), estimated from above by
$C_{\alpha,A,m,\chi,\psi}\|u^{\delta} - u\|_{L^2}$.
The first term is equal to
\begin{align*}%
&\|\chi \big([(H_{0,m})^{\alpha},\psi] + \psi (H_{0,m})^{\alpha}\big)
              (u^{\delta} - u)\|_{L^1}\\
&\leq \|\chi [(H_{0,m})^{\alpha},\psi](u^{\delta} - u)\|_{L^1}
 + \|\chi \psi
  [(H_{0,m})^{\alpha} u^{\delta} - (H_{0,m})^{\alpha}u]\|_{L^1}\\
&\leq C_{\alpha,0,m,\chi,\psi}\|u^{\delta} - u\|_{L^2}
  + \|((H_{0,m})^{\alpha}u)^{\delta} - (H_{0,m})^{\alpha}u\|_{L^1},
\end{align*}%
where we have used for the first term Lemma 2.5 (2.19) for $A=0$
and for the second the fact
that $(H_{0,m})^{\alpha} u^{\delta} = ((H_{0,m})^{\alpha} u)^{\delta}$
since by assumption, $(H_{0,m})^{\alpha}u$ is locally in
$L^1$ and $u \in L^2$. It follows that
\begin{align*}%
&\|(H_{A,m})^{\alpha} u^{\delta} - (H_{A,m})^{\alpha} u\|_{L^1(K)}\\
&\leq C_{\alpha,0,m,\chi,\psi}\|u^{\delta} - u\|_{L^2}
 +\|\psi\|_{L^{\infty}}\|((H_{0,m})^{\alpha}u)^{\delta}
                                  - (H_{0,m})^{\alpha}u\|_{L^1}
 + C_{\alpha,A,m,\chi,\psi}\|u^{\delta} - u\|_{L^2},
\end{align*}%
which approaches zero as $\delta \downarrow 0$.
This proves Lemma 3.1.
\qed

\begin{lemma}
Let $0<\alpha \leq 1$. Let $u \in L^2({\Bbb R}^d)$ and
$H_{A,m} u \in L^1_{\text{\rm loc}}({\Bbb R}^d)$. Then
 $(H_{A,m})^{\alpha}u = [(-i\nabla -A)^2+m^2]^{\alpha/2}u$
is also in $L^1_{\text{\rm loc}}({\Bbb R}^d)$, and $\{ (H_{A,m})^{\alpha}u \}$
 converges to $H_{A,m} u$ in $L^1_{\text{\rm loc}}({\Bbb R}^d)$
as $\alpha \uparrow 1$. Namely,
for any $\psi \in C_0^{\infty}({\Bbb R}^d)$,
$\|\psi(H_{A,m})^{\alpha}u\|_{L^1}$ is uniformly bounded for $0<\alpha \leq 1$,
and $\{ \psi(H_{A,m})^{\alpha}u \}$ converges to $\psi  H_{A,m} u$
in $L^1({\Bbb R}^d)$ as $\alpha \uparrow 1$.
\end{lemma}

\bigskip
{\it Proof}.
Let $0<\alpha <1$.
To begin with, suppose with $\psi \in C_0^{\infty}({\Bbb R}^d)$
that some $u \in L^2({\Bbb R}^d)$ satisfies the equality
\begin{equation}%
\psi (H_{A,m})^{\alpha}u
= (H_{A,m})^{-(1-\alpha)}\psi H_{A,m} u
        + [\psi, (H_{A,m})^{-(1-\alpha)}]H_{A,m} u.
\end{equation}
This holds at least if $u \in D[H_{A,m}]$, and hence, in particular, if
$u = \phi \in C_0^{\infty}({\Bbb R}^d)$. Note here that
$(H_{A,m})^{\alpha}$ has $D[H_{A,m}]$ as an operator core, while
$H_{A,m}$ has $C_0^{\infty}({\Bbb R}^d)$ as an operator core.

Now, let $u \in L^2({\Bbb R}^d)$ with $H_{A,m} u
\in L^1_{\text{\rm loc}}({\Bbb R}^d)$, just what is assumed by
Lemma 3.2. The first term on the right-hand side of (3.2)
is in $L^1({\Bbb R}^d)$, since by Lemma 2.1 (i) with $p=1$,
$(H_{A,m})^{-(1-\alpha)}$ is a bounded operator which is a contraction
mapping $L^1({\Bbb R}^d)$ into $ L^1({\Bbb R}^d)$,
 bounded uniformly for $0<\alpha\leq 1$ and strongly
continuous there, so long as $m>0$.
 The term on the left-hand side of (3.2) exists as a distribution.
The second term on the right-hand side lies in the dual space of the space
$D[H_{A,m}]$, considered as a Hilbert space equipped with the graph norm
$\|v\|_{L^2}^2 + \|H_{A,m}v\|^2$. Here recall (2.1) and note that
for $\phi \in C_0^{\infty}({\Bbb R}^d)$
$$\|(H_{A,m})^{\alpha}\phi\|_{L^2}
= \|(H_{A,m})^{-(1-\alpha)}\,H_{A,m}\phi\|_{L^2}
\leq \|H_{A,m}\phi\|_{L^2}.
$$
Thus all the three terms on the
left- and right-hand sides of (3.2) exist also as distributions.

To show the assertion of the lemma, take a $C^{\infty}$ cutoff function
$\chi$ with compact support,
a similar one of which has already been used before, such that
$0\leq \chi(x) \leq 1$ in ${\Bbb R}^d$ with
$\chi(x)=1$ on $\text{\rm supp}\, \psi$.
As $\psi = \chi\psi$ holds,
so does $\psi (H_{A,m})^{\alpha}u = \chi\,\psi (H_{A,m})^{\alpha}u$.
Then consider the (3.2) multiplied by $\chi$, i.e.
\begin{equation}%
\psi (H_{A,m})^{\alpha}u
= \chi\, (H_{A,m})^{-(1-\alpha)}\psi H_{A,m} u
         + \chi\, [\psi, (H_{A,m})^{-(1-\alpha)}]H_{A,m} u.
\end{equation}
The first term on the right of (3.2)
(and hence (3.3)) converges to $\psi H_{A,m}u$
as $\alpha \uparrow 1$, since $(H_{A,m})^{-(1-\alpha)}$
is an operator on $ L^1({\Bbb R}^d)$,
 bounded uniformly for $0<\alpha\leq 1$ and strongly continuous there,
so long as $m>0$.
So we have only to show the second term of (3.3),
i.e. $\chi [\psi, (H_{A,m})^{-(1-\alpha)}]H_{A,m} u$ lie
in $L^1({\Bbb R}^d)$, being uniformly bounded, and converge to $0$
in $L^1$ as $\alpha \uparrow 1$.

Use formula (2.2) to rewrite this second term on the right of (3.3) as
\begin{align}%
&\chi [\psi, (H_{A,m})^{-(1-\alpha)}]H_{A,m} u      \nonumber\\
&= \frac1{\Gamma(\frac{1-\alpha}2)}\int_0^{\infty} dt\,
   t^{\frac{1-\alpha}2 -1}\,
  \chi \big[\psi e^{-t(H_{A,m})^2} - e^{-t(H_{A,m})^2}\psi\big]H_{A,m} u
                                                 \nonumber\\
&= -\frac1{\Gamma(\frac{1-\alpha}2)}\int_0^{\infty} dt\,
  t^{-\frac{1+\alpha}2}\, \chi \int_0^1 d\theta\, \frac{d}{d\theta}
  \Big[e^{-\theta t(H_{A,m})^2}\psi e^{-(1-\theta)t(H_{A,m})^2}\Big]H_{A,m} u
                                                  \nonumber\\
&= \frac1{\Gamma(\frac{1-\alpha}2)}\int_0^{\infty} dt\,
 t^{\frac{1-\alpha}2}
 \big[(H_{A,m})^2,\psi \big]e^{-(1-\theta)t(H_{A,m})^2}H_{A,m} u.
\end{align}%
Recall identity (2.10) for the commutator $\big[(H_{A,m})^2,\psi \big]$,
indeed, the first of the two expressions there and
substitute it into the $\big[(H_{A,m})^2,\psi \big]$ in the last member
of (3.4). Then
\begin{align}%
&\chi [\psi, (H_{A,m})^{-(1-\alpha)}]H_{A,m} u \nonumber\\
&= \frac1{\Gamma(\frac{1-\alpha}2)}\int_0^{\infty} dt\,
    t^{\frac{1-\alpha}2}\,\int_0^1 d\theta\,
    \chi\, e^{-\theta t(H_{A,m})^2}\, (\Delta\psi)\,
    e^{-(1-\theta)t(H_{A,m})^2}H_{A,m} u  \nonumber\\
 &\quad +\frac2{\Gamma(\frac{1-\alpha}2)}\int_0^{\infty} dt\,
    t^{\frac{1-\alpha}2}\,\int_0^1 d\theta\,
    \chi\, e^{-\theta t(H_{A,m})^2} (i\nabla +A)(i\nabla\psi)
    e^{-(1-\theta)t(H_{A,m})^2}H_{A,m} u \nonumber\\
&=: I_5u + I_6u.                   
\end{align}%
We estimate the $L^1$ norms of $I_5u$ and $I_6u$ in (3.5).

First for $I_5u$, integrate its absolute value in $x$, we have
by the Schwarz inequality
\begin{align*}%
\|I_5u\|_{L^1}
&\leq \frac1{\Gamma(\frac{1-\alpha}2)}\int_0^{\infty} dt\,
    t^{\frac{1-\alpha}2}\,\int_0^1 d\theta\,
    \big\|\chi\, e^{-\theta t(H_{A,m})^2}\, (\Delta\psi)\,
    e^{-(1-\theta)t(H_{A,m})^2}H_{A,m} u\big\|_{L^1}\\
&\leq \frac1{\Gamma(\frac{1- \alpha}2)}\int_0^{\infty} dt\,
    t^{\frac{1-\alpha}2}\,\int_0^1 d\theta\,
    \big\|\chi\big\|_{L^2}\\
 &\qquad\qquad\qquad\qquad \times
    \big\|e^{-\theta t(H_{A,m})^2}\, (\Delta\psi)\,
    e^{-(1-\theta)t(H_{A,m})^2}H_{A,m} u\big\|_{L^2}\\
&\leq \frac1{\Gamma(\frac{1-\alpha}2)}\int_0^{\infty} dt\,
    t^{\frac{1-\alpha}2}\,\int_0^1 d\theta\,
    \big\|\chi\big\|_{L^2}
    \big\|e^{-\theta t(H_{A,m})^2}\big\|_{L^2\rightarrow L^2}\,\\
 &\qquad\qquad\qquad\qquad \times
    \big\|\Delta\psi\big\|_{L^{\infty}}\,
    \big\|e^{-(1-\theta)t(H_{A,m})^2}H_{A,m}
    \big\|_{L^2\rightarrow L^2}\big\|u\big\|_{L^2} \\
&\leq \frac1{\Gamma(\frac{1-\alpha}2)}\int_0^{\infty} dt\,
    t^{\frac{1-\alpha}2}\,\int_0^1 d\theta\,
    \big\|\chi\big\|_{L^2} e^{-\frac{m^2}2 \theta t}
    \big\|e^{-\frac{\theta t}2(H_{A,m})^2}\big\|_{L^2\rightarrow L^2}\,\\
 &\qquad\qquad\qquad\qquad \times
    \big\|\Delta\psi\big\|_{L^{\infty}}\, e^{-\frac{m^2}2 (1-\theta) t}
    \big\|e^{-\frac{(1-\theta)t}2(H_{A,m})^2}H_{A,m}
    \big\|_{L^2\rightarrow L^2}\big\|u\big\|_{L^2}\,.
\end{align*}%
Then by Lemma 2.1 (iii) we have the bound
\begin{eqnarray}%
\|I_5u\|_{L^1}
&\leq& \frac1{\Gamma(\frac{1-\alpha}2)}\int_0^{\infty} dt\,
    t^{\frac{1-\alpha}2} e^{-\tfrac{m^2}2 t}\,\int_0^1 d\theta\,
    \big\|\chi\big\|_{L^2}\|\Delta\psi\big\|_{L^{\infty}}\,
    \big(\tfrac{1}{2e(\tfrac{1-\theta}2)t}\big)^{1/2}\big\|u\big\|_{L^2}
                                                    \nonumber\\
&\leq& \frac1{\Gamma(\frac{1-\alpha}2)}\int_0^{\infty} dt\,
    t^{-\frac{\alpha}2}e^{-\tfrac{m^2}2 t}\,
    (2e)^{-1/2}\int_0^1 \tfrac{d\theta}{(\frac{1-\theta}2)^{1/2}}\,
    \big\|\chi\big\|_{L^2}\|\Delta\psi\big\|_{L^{\infty}}
    \big\|u\big\|_{L^2}                             \nonumber\\
&\leq& 2^{\frac32}(2e)^{-1/2}
   \frac{\Gamma\big(\frac{2-\alpha}2\big)}{\Gamma\big(\frac{1-\alpha}2\big)}
   \big(\tfrac2{m^2}\big)^{\frac{2-\alpha}2}
   \big\|\chi\big\|_{L^{2}} \big\|\Delta\psi\big\|_{L^{\infty}}
   \big\|u\big\|_{L^2}.
\end{eqnarray}

Next for $I_6 u$,  we are going to show a similar bound
\begin{equation}%
\|I_6u\|_{L^1}
\leq 4\pi(2e)^{-1/2}
  \frac{\Gamma\big(\frac{2-\alpha}2\big)}{\Gamma\big(\frac{1-\alpha}2\big)}
  \big(\tfrac{2}{m^2}\big)^{\tfrac{2-\alpha}2}
  C_{\chi,A} \big\|\nabla\psi\big\|_{L^{\infty}}\big\|u\big\|_{L^2},
\end{equation}
with a constant, depending only on $\chi$ and $A$,
\begin{equation}%
C_{\chi,A} := \big[\|\nabla \chi\|_{L^2}^2
             +m^2\|\chi\|_{L^2}^2 + \|\chi\, A\|_{L^2}^2\big]^{1/2},
\end{equation}
 which is bounded since
$A \in L^2_{\text{\rm loc}}({\Bbb R}^d)$.
The proof is to integrate the absolute value of $I_{6}u$ in $x$ to get
\begin{align}%
\|I_{6}u\|_{L^1}
&\leq
\frac2{\Gamma(\frac{1-\alpha}2)}\int_0^{\infty} dt\,
  t^{\frac{1-\alpha}2}\,\int_0^1 d\theta\,
  X_{A,m}(t, \theta; \chi, \psi, u),      \\ 
\intertext{\rm where we put}
X_{A,m}(t, \theta; \chi, \psi, u)
&= \big\|\chi\, e^{-\theta t(H_{A,m})^2}(i\nabla+A) (i\nabla\psi)
  e^{-(1-\theta)t(H_{A,m})^2} H_{A,m} u\big\|_{L^1}\,.
\end{align}
Somewhat crucial is the estimate of
$X_{A,m}(t, \theta; \chi, \psi, u)$
in (3.10) which we are going to do,
where the parentheses $\big( \cdot, \cdot\big)$ below
stand for the $L^2$ inner product:
\begin{align}
&X_{A,m}(t, \theta; \chi, \psi, u)                      \nonumber\\
&= \big|\big(\chi, e^{-\theta t(H_{A,m})^2}(i\nabla+A) (i\nabla\psi)
   e^{-(1-\theta)t(H_{A,m})^2} H_{A,m} u\big)\big|      \nonumber\\
&= \big|\big(\chi, e^{-\theta t(H_{A,m})^2} H_{A,m}\cdot
       (H_{A,m})^{-1}(i\nabla+A) (i\nabla\psi)
        e^{-(1-\theta)t(H_{A,m})^2} H_{A,m} u\big)\big| \nonumber\\
&= \big|\big(e^{-\theta t(H_{A,m})^2} H_{A,m}\, \chi,
        (H_{A,m})^{-1}(i\nabla+A) (i\nabla\psi)
        e^{-(1-\theta)t(H_{A,m})^2} H_{A,m} u\big)\big| \nonumber\\
&\leq \big\|e^{-\theta t(H_{A,m})^2} H_{A,m}\, \chi\|_{L^2}\,
      \big\|(H_{A,m})^{-1}(i\nabla+A) (i\nabla\psi)
          e^{-(1-\theta)t(H_{A,m})^2} H_{A,m} u\big\|_{L^2}.
\end{align}
In the last member of (3.11), the first factor and the second are
estimated as follows:
\begin{align}
\big\|e^{-\theta t(H_{A,m})^2} H_{A,m}\, \chi\|_{L^2}
&\leq e^{-m^2 \theta t} \big\|H_{A,m} \chi\big\|_{L^2}             \nonumber\\
&= e^{-m^2 \theta t}\big[\sum_{j=1}^d
 \|(i\partial_j +A_j)\chi\|_{L^2}^2 +m^2\|\chi\|_{L^2}^2\big]^{1/2}\nonumber\\
&\leq e^{-m^2 \theta t}\big[\|\nabla \chi\|_{L^2}^2
             +\|\chi\, A\|_{L^2}^2+m^2\|\chi\|_{L^2}^2\big]^{1/2}  \nonumber\\
&= e^{-m^2 \theta t}\, C_{\chi,A}\,;
\end{align}
\begin{align}%
&\big\|(H_{A,m})^{-1}(i\nabla+A) (i\nabla\psi)
          e^{-(1-\theta)t(H_{A,m})^2} H_{A,m} u\big\|_{L^2 } \nonumber\\
&\leq \big\|(H_{A,m})^{-1}(i\nabla+A)\big\|_{L^2\rightarrow L^2}
      \big\|(i\nabla\psi)\big\|_{L^{\infty}}
      \big\|e^{-(1-\theta)t(H_{A,m})^2} H_{A,m}
                                     \big\|_{L^2\rightarrow L^2}
      \big\|u\big\|_{L^2}                                    \nonumber\\
&\leq \big\|\nabla\psi\big\|_{L^{\infty}} e^{-\frac{m^2}2 (1-\theta) t}
      \big\|e^{-\frac{(1-\theta)t}2(H_{A,m})^2} H_{A,m}
               \big\|_{L^2\rightarrow L^2}\big\|u\big\|_{L^2} \nonumber\\
&\leq \big\|\nabla\psi\big\|_{L^{\infty}} e^{-\frac{m^2}2 (1-\theta) t}
      \big(\tfrac{1}{2e \frac{(1-\theta) t}2}\big)^{1/2}
                                                   \big\|u\big\|_{L^2}.
\end{align}
Here in (3.12) and (3.13) we have used (2.1), Lemma 2.1 (iii),
and the estimate $\big\|(H_{A,m})^{-1}(i\nabla+A)\big\|_{L^2\rightarrow L^2} \leq 1$.
From (3.12) and (3.13) we obtain
\begin{align*}%
\|I_6u\|_{L^1}
&\leq  \frac2{\Gamma(\frac{1-\alpha}2)}\int_0^{\infty} dt\,
  t^{\frac{1-\alpha}2}e^{-\frac{m^2}2 t}\,\int_0^1 d\theta\,
  \big(\tfrac{1}{2e \frac{(1-\theta) t}2}\big)^{1/2}
  C_{\chi,A}\big\|\nabla\psi\big\|_{L^{\infty}}\big\|u\big\|_{L^2}
                                                \nonumber\\
&\leq \frac2{\Gamma(\frac{1-\alpha}2)}\int_0^{\infty} dt\,
  t^{-\frac{\alpha}2}e^{-\frac{m^2}2t}\,
  \int_0^1 \tfrac{d\theta}{(1-\theta)^{1/2}}\,
  (2e)^{-1/2} C_{\chi,A}\big\|\nabla\psi\big\|_{L^{\infty}}
  \big\|u\big\|_{L^2}\,,      \nonumber\\
\end{align*}
of which the last member yields (3.7) with (3.8).

Thus, taking (3.5) into account and putting together (3.6) and (3.7),
we see the $L^1$ norm of the second term on the right-hand side of
(3.3) is estimated as
\begin{align}%
&\|\chi[\psi, (H_{A,m})^{-(1-\alpha)}]H_{A,m} u\|_{L^1}\nonumber\\
&\leq \|I_5 u\|_{L^1} + \|I_{6} u\|_{L^1}              \nonumber\\
&\leq \frac{\Gamma\big(\frac{2-\alpha}2\big)}{\Gamma\big(\frac{1-\alpha}2\big)}
   \big(\tfrac2{m^2}\big)^{\frac{2-\alpha}2}
   \Big[2^{\frac32}(2e)^{-1/2}
   \big\|\chi\big\|_{L^2} \big\|\Delta\psi\big\|_{L^2}
   +  4\pi (2e)^{-1} C_{\chi,A} \|\nabla\psi\|_{L^{\infty}} \Big]                                                                      \|u\|_{L^2}.
\end{align}
Since the last member of (3.14) tends to zero as $\alpha \uparrow 1$,
 because  $\Gamma(z) \uparrow \infty$ as $z \downarrow 0$ and hence
 $\frac1{\Gamma(\frac{1-\alpha}2)} \rightarrow 0$ as
$\alpha \uparrow 1$,
we see the left-hand side is uniformly bounded for $0<\alpha <1$,
and convergent to zero as $\alpha \uparrow 1$. This shows
the desired assertion of Lemma 3.2. \qed

\bigskip

Now we are in a position to prove Theorem 1.1.

\bigskip
{\it Completion of Proof of Theorem 1.1}.

\medskip
As (1.2) and (1.3) are equivalent, we have only to show (1.3).
The proof is divided into three parts,
(i) the case where $m>0$ and $0<\alpha<1$,
(ii) the case where $m>0$ and $\alpha=1$,
(iii) the case where $m=0$ and $\alpha=1$.

\medskip
(i) {\it The case where $m>0$ and $0<\alpha<1$}.
We prove in two steps, treating first the step (i-I) for
$u \in (C^{\infty}\cap L^2)({\Bbb R}^d)$, and next the step (i-II)
for general $u \in L^2$ with $(H_{A,m})^{\alpha}u \in L^1_{\text{\rm loc}}$.

\bigskip
(i-I) {\it For $u \in (C^{\infty}\cap L^2)({\Bbb R}^d) \, (0<\alpha<1)$}.

For a function $v(x) \in C^{\infty}({\Bbb R}^d)$ and
$\varepsilon>0$, put
$v_{\varepsilon}(x)= \sqrt{|v(x)|^2+\varepsilon^2}$. Then note that
$v_{\varepsilon}(x) \geq \varepsilon$, and, since
$v_{\varepsilon}(x)^2= |v(x)|^2+\varepsilon^2$,
we have
{\begin{equation}
 -|v(x)||v(x+y)|+|v(x)|^2
 \geq -v_{\varepsilon}(x)v_{\varepsilon}(x+y)+v_{\varepsilon}(x)^2.
\end{equation}
}
Then we will show that
$u_{\varepsilon} = \sqrt{|u|^2+\varepsilon^2}, \,\,\varepsilon >0$,
satisfies that
\begin{equation}%
 \hbox{\rm Re}[\overline{u(x)}([(H_{A,m})^{\alpha}-m^{\alpha}] u)(x)]
              \geq u_{\varepsilon}(x)
                   \big([(H_{0,m})^{\alpha}-m^{\alpha}]
                   (u_{\varepsilon}-\varepsilon)\big)(x),
 \quad \hbox{\rm  pointwise} \,\, \hbox{\rm  a.e.},
\end{equation}
which amounts to the same thing as
\begin{equation}%
 \hbox{\rm Re}\Biggl[\frac{\overline{u(x)}}{u_{\varepsilon}(x)}
      ([(H_{A,m})^{\alpha}-m^{\alpha}] u)(x)\Biggr]
      \geq  \big([(H_{0,m})^{\alpha}-m^{\alpha}]
                 (u_{\varepsilon}-\varepsilon)\,\big)(x),
\end{equation}
pointwise a.e., and thus in distributional sense. Here note that the function
$u_{\varepsilon}-\varepsilon$ is nonnegative, $C^{\infty}$ and has the same
compact support as $u$.

\medskip
We show (3.16) or
{(3.17)} first for $u \in C_0^{\infty}({\Bbb R}^d)$ and then for 
$u \in (C^{\infty}\cap L^2)({\Bbb R}^d)$. To do so, we employ analogous
arguments as used in [I93, p.223, Lemma 2] for the case $\alpha=1$, i.e.
for $H_{A,m}-m$. We will use the same notation $S$ as in Section 2 for
the selfadjoint operator
$(-i\nabla-A(x))^2+m^2$ in $L^2({\Bbb R}^d)$, which may be
considered as the magnetic nonrelativistic Schr\"odinger operator
with mass $\frac12$ with constant scalar potential $m^2$.
Then we have $H_{A,m} = S^{\frac12}$. Since the domain of
$H_{A,m}$ includes $C_0^{\infty}({\Bbb R}^d)$ as a operator core, 
the operator $[H_{A,m}-m]u$ can be written as
$\text{s-}\lim_{t \downarrow 0}t^{-1}\big(1-e^{-t[H_{A,m}-m]}\big)u$.
It is known from the theory of fractional powers of a linear operator
(e.g. see [Y78, IX, 11, pp.259--261]) that
the semigroup $e^{-t[(H_{A,m})^{\alpha}-m^{\alpha}]}$ with generator
$(H_{A,m})^{\alpha} = S^{\frac{\alpha}2}$
is obtained from the semigroup $e^{-tS}$ with generator $S$ as
\begin{equation}
 e^{-t[(H_{A,m})^{\alpha}-m^{\alpha}]}u =
\left\{
 \begin{array}{ll}
     e^{m^{\alpha}t}\int_0^{\infty}
                   f_{t,\frac{\alpha}2}(\lambda)e^{-\lambda S}u
                   d\lambda, &\quad t>0,\\
       u, &\quad t=0,
\end{array}\right.
\end{equation}
where for $t>0$ and $0<\alpha \leq 1$, $f_{t,\frac{\alpha}2}(\lambda)$ is
a nonnegative function of exponential growth in $\lambda \in {\Bbb R}$
given by
\begin{equation}
f_{t,\frac{\alpha}2}(\lambda)=
\left\{
\begin{array}{ll}
     (2\pi i)^{-1}\int_{\sigma-i\infty}^{\sigma+i\infty}
       e^{z\lambda-tz^{\frac{\alpha}2}}dz, &\quad \lambda \geq 0,\\
                             0, &\quad \lambda <0,
\end{array}\right.
\end{equation}
with $\sigma>0$, where
the branch of $z^{\frac{\alpha}2}$ is so taken that
$\text{Re}\, z^{\frac{\alpha}2} > 0$ for $\text{Re}\, z > 0$.
In passing, we note that equation (3.18) is valid even for $1<\alpha < 2$,
though we don not need this case in the present paper.

We continue our preceding arguments and recall that
$
 |e^{-tS}u| \leq e^{-t[-\Delta+m^2]}|u|
$
pointwise a.e., what is referred to in (2.17).
It follows with (3.18), (3.19), that
\begin{align}
 |e^{-t[(H_{A,m})^{\alpha}-m^{\alpha}]}u|
&\leq e^{m^{\alpha}t}\int_0^{\infty}
     f_{t,\frac{\alpha}2}(\lambda)|e^{-\lambda S}u|d\lambda \nonumber\\
&\leq e^{m^{\alpha}t}\int_0^{\infty}
             f_{t,\frac{\alpha}2}(\lambda)e^{-\lambda (-\Delta+m^2)}
                |u|d\lambda                                 \nonumber\\
&= e^{-t[(H_{0,m})^{\alpha}-m^{\alpha}]}|u|,  
\end{align}
poitwise a.e.
Hence for $t>0$
\begin{equation}
 \text{Re}\Big[\,\overline{u(x)}
   \Big(\frac{1-e^{-t[(H_{A,m})^{\alpha}-m^{\alpha}]}}{t} u\Big)(x)\Big]
\geq |u(x)| \Big(\frac{1-e^{-t[(H_{0,m})^{\alpha}-m^{\alpha}]}}{t}|u|\Big)(x),
\end{equation}
poitwise a.e.
Now put $n^{m,\alpha}(t,y) := \frac1{t}k_0^{m,\alpha}(t,y)$, taking account of
the relation (2.7) between the integral kernel
$k_0^{m,\alpha}(t,y)$ of $e^{-t[(H_{0,m})^{\alpha}-m^{\alpha}]}$ and
the density (function) $n^{m,\alpha}(y)$ of the L\'evy measure.

Then we see, by (2.6), that the right-hand side of (3.21) equal to
\begin{eqnarray*}%
&& |u(x)|\int_{|y|>0} [|u(x)|-|u(x+y)|] \,\frac{k_0^{m,\alpha}(t,y)}{t} dy\\
&\quad& = -\int_{|y|>0}
        \big[\,|u(x)||u(x+y)|-|u(x)|^2 \, \big] \,n^{m,\alpha}(t,y) dy\\
&\quad& \geq
-\int_{|y|>0} \big[\,u_{\varepsilon}(x)
   u_{\varepsilon}(x+y)-u_{\varepsilon}(x)^2\,\big] \,n^{m,\alpha}(t,y) dy\\
&\quad&= u_{\varepsilon}(x)\Big[-
  \int_{|y|>0}  \big[\,u_{\varepsilon}(x+y)-u_{\varepsilon}(x)
  -I_{\{|y|<1\}}y\cdot\nabla u_{\varepsilon}(x) \,\big] \,
   n^{m,\alpha}(t,y) dy\Big],
\end{eqnarray*}%
for every $\varepsilon >0$,
where we used (3.15) and the $y$-rotational invariance of
$k_0^{m,\alpha}(t,y)$ or $n^{m,\alpha}(t,y)$. Notice
the integral $\Big[-\int_{|y|>0} \cdots\Big]$ of the last member is
equal to that with $(u_{\varepsilon} -\varepsilon)$ in place of
$u_{\varepsilon}$, i.e.
$$
\Big(\frac{1-e^{-t[(H_{0,m})^{\alpha}-m^{\alpha}]}}{t}
     (u_{\varepsilon} -\varepsilon)\Big)(x).
$$
Thus we have from (3.21)
\begin{equation}%
\text{Re}\, \Big[\overline{u(x)}
    \Big(\frac{1-e^{-t[(H_{A,m})^{\alpha}-m^{\alpha}]}}{t}u\Big)(x)\Big]
\geq u_{\varepsilon}(x)\Big(\frac{1-e^{-t[(H_{0,m})^{\alpha}-m^{\alpha}]}}{t}
     (u_{\varepsilon} -\varepsilon)\Big)(x).
\end{equation}
Then letting $t \downarrow 0$ on both sides of (3.22), we obtain (3.16).
Indeed, recalling the function $u_{\varepsilon} -\varepsilon$ has
compact support,
the right-hand side tends to that of (3.16).
For the left-hand side, since $u$ is in the domain of
$(H_{A,m})^{\alpha}-m^{\alpha}$, we have
$t^{-1}[1-e^{-t[(H_{A,m})^{\alpha}-m^{\alpha}]}]u
        \rightarrow [(H_{A,m})^{\alpha}-m^{\alpha}]u$ in $L^2$,
and pointwise a.e. by passing to a subsequence.
This shows (3.16)/(3.17) for $u \in C_0^{\infty}({\Bbb R}^d)$.

\bigskip
Next we show (3.16)/(3.17) when $u \in (C^{\infty} \cap L^2)({\Bbb R}^d)$.
Take a sequence $\{u_n\} \in  C_0^{\infty}({\Bbb R}^d)$ such that
$u_n \rightarrow u$ in $(C^{\infty} \cap L^2)({\Bbb R}^d)$,
i.e. in the topology of $C^{\infty}({\Bbb R}^d)$ as well as
in the norm of $L^2({\Bbb R}^d)$, as $n\rightarrow \infty$.
Then from the case $u \in C_0^{\infty}({\Bbb R}^d)$ above,
we have for all $\varepsilon >0$
$$
 \hbox{\rm Re}
\Biggl[\frac{\overline{u_n(x)}}
   {u_{n,\varepsilon}(x)}\big([(H_{A,m})^{\alpha}-m^{\alpha}] u_n\big)(x)\Biggr]
\geq  \big([(H_{0,m})^{\alpha}-m^{\alpha}] (u_{n,\varepsilon}-\varepsilon)\big)(x),
$$
pointwise, and hence
for any $\psi \in C_0^{\infty}({\Bbb R}^d)$ with $\psi(x) \geq 0$,
$$
 \hbox{\rm Re}
 \Big\langle \psi, \,\frac{\overline{u_n}}
 {u_{n,\varepsilon}}\big([(H_{A,m})^{\alpha}-m^{\alpha}] u_n\big) \Big\rangle
 \geq \Big\langle \psi, \,
  [(H_{0,m})^{\alpha}-m^{\alpha}] (u_{n,\varepsilon}-\varepsilon) \Big\rangle
$$
for all $\varepsilon >0$. Here the bilinear inner product
$\langle\cdot, \cdot\rangle$ is an integral with respect to the Lebesgue measure $dx$,
and also considered as the bilinear inner product
between the dual pair of the test functions and the distributions:
$\big\langle  C_0^{\infty}({\Bbb R}^d), \, {\cal D}'({\Bbb R}^d) \big\rangle$.
Therefore
$$
 \hbox{\rm Re}
 \Big\langle [(H_{A,m})^{\alpha}-m^{\alpha}]
 \Big(\frac{\overline{u_n}}{u_{n,\varepsilon}} \psi\Big), \,
 u_n \Big\rangle
 \geq \Big\langle [(H_{0,m})^{\alpha}-m^{\alpha}] \psi,\,
 u_{n,\varepsilon}-\varepsilon \Big\rangle.
$$
Since we have that $u_n \rightarrow u$  and
$u_{n,\varepsilon} \rightarrow u_{\varepsilon}$ in
$(C^{\infty}\cap L^2)({\Bbb R}^d)$ as $n\rightarrow \infty$, we have that
$\big(\tfrac{\overline{u_n}}{u_{n,\varepsilon}}\big)\psi
\rightarrow
\big(\tfrac{\overline{u}}{u_{\varepsilon}}\big)\psi
$. It follows by Lemma 2.2 that
$$
[(H_{A,m})^{\alpha}-m^{\alpha}]
\Big(\frac{\overline{u_n}}{u_{n,\varepsilon}}\Big)\psi
\rightarrow
([H_{A,m})^{\alpha}-m^{\alpha}]
\Big(\frac{\overline{u}}{u_{\varepsilon}}\Big)\psi
$$
in $L^2$ as $n\rightarrow \infty$, so that
$$
 \hbox{\rm Re}
 \Big\langle [(H_{A,m})^{\alpha}-m^{\alpha}]
      \Big(\frac{\overline{u}}{u_{\varepsilon}}\psi\Big), \, u \Big\rangle
     \geq \Big\langle [(H_{0,m})^{\alpha}-m^{\alpha}] \psi, \,
   (u_{\varepsilon}-\varepsilon) \Big\rangle.
$$
Thus we obtain
\begin{equation}%
 \hbox{\rm Re}
\Biggl[\frac{\overline{u(x)}}
{u_{\varepsilon}(x)}\big([(H_{A,m})^{\alpha}-m^{\alpha}] u\big)(x)\Biggr]
\geq  \big([(H_{0,m})^{\alpha}-m^{\alpha}] (u_{\varepsilon}-\varepsilon)\big)(x),
\end{equation}
pointwise a.e., and thus in distributional sense, and hence (3.17) follows 
for $u \in (C^{\infty} \cap L^2)({\Bbb R}^d)$.

\medskip
(i-II) {\it For general $u \in L^2({\Bbb R}^d)$ with
$(H_{A,m})^{\alpha} u \in L^1_{\text{\rm \rm loc}}({\Bbb R}^d)\, (0<\alpha<1)$}.

Put $u^{\delta} = \rho_{\delta} *u$.
Then $u^{\delta} \in C^{\infty} \cap L^2$, so by (3.23) in step (i-I) above
\begin{equation}%
 \hbox{\rm Re}
 \Biggl[\frac{\overline{u^{\delta}}}{(u^{\delta})_{\varepsilon}}
 \big([(H_{A,m})^{\alpha}-m^{\alpha}] u^{\delta}\big)\Biggr]
 \geq  [(H_{0,m})^{\alpha}-m^{\alpha}]
       \big((u^{\delta})_{\varepsilon}-\varepsilon\big),
\end{equation}
pointwise a.e., and also in distributional sense,
for all $\varepsilon >0$ and all $\delta>0$.

We first, for fixed $\varepsilon >0$, let $\delta \downarrow 0$, and next
$\varepsilon \downarrow 0$. In fact, if $\delta \downarrow 0$, then
$u^{\delta} \rightarrow u$ in $L^2$ as well as a.e. by passing to a subsequence
of $\{ u^{\delta} \}$. Hence
$\overline{u^{\delta}}/(u^{\delta})_{\varepsilon} \rightarrow
 \overline{u}/u_{\varepsilon}$ a.e. and by Lemma 3.1,
$(H_{A,m})^{\alpha} u^{\delta} \rightarrow (H_{A,m})^{\alpha} u$ locally in $L^1$,
{ and therefore also a.e. by passing to a subsequence.
}
Since
$\big|\frac{\overline{u^{\delta}}}{(u^{\delta})_{\varepsilon}}\big| \leq 1$,
it follows by the Lebesque dominated convergence theorem that
on the left-hand side of (3.24),
$$
 \frac{\overline{u^{\delta}}}{(u^{\delta})_{\varepsilon}}
 [(H_{A,m})^{\alpha}-m^{\alpha}]u^{\delta} \rightarrow
 \frac{\overline{u}}{u_{\varepsilon}}  [(H_{A,m})^{\alpha}-m^{\alpha}]u
$$
locally in $L^1$ as $\delta \downarrow 0$. On the other hand,
for the right-hand side, since
$$
\big|\big((u^{\delta})_{\varepsilon} -\varepsilon\big)
           - (u_{\varepsilon} -\varepsilon)\big|
\leq \big|(u^{\delta})_{\varepsilon} - u_{\varepsilon}\big|
\leq  \big||u^{\delta}| - |u|\big|\leq |u^{\delta} - u|,
$$
we have
$(H_{0,m})^{\alpha} \big((u^{\delta})_{\varepsilon}-\varepsilon\big)
   \rightarrow
    (H_{0,m})^{\alpha} (u_{\varepsilon}-\varepsilon)
$
in ${\cal D}'$ (in distributional sense). This shows that (3.23)
holds for $u \in L^2({\Bbb R}^d)$ with
$(H_{A,m})^{\alpha} u \in L^1_{\text{\rm loc}}({\Bbb R}^d)$.
Next let $\varepsilon \downarrow 0$. Then
$\overline{u}/u_{\varepsilon} \rightarrow \text{sgn}\, u$ a.e. with
$|\overline{u}/u_{\varepsilon}|\leq 1$, so that the left-hand side
of (3.23) converges to
$\text{Re} ((\text{sgn}\, u)[H_{A,m} -m]u)$ a.e.,
while the right-hand side of (3.23)
converges to $[(H_{0,m})^{\alpha}-m^{\alpha}]|u|$ in ${\cal D}'$.
Thus we get (3.1), showing the desired inequality for $0<\alpha<1$.

\medskip
(ii) {\it The case where $m>0$ and $\alpha=1$}.

Once the inequality (3.1) is established for $0<\alpha<1$,
we let $\alpha \uparrow 1$, with $u \in L^2({\Bbb R}^d)$ with
$H_{A,m} u \in L^1_{\text{\rm loc}}({\Bbb R}^d)$.
Then, as $\alpha \uparrow 1$, by Lemma 3.2 we have
$(H_{A,m})^{\alpha}u \rightarrow H_{A,m} u$ in $L^1_{\text{\rm loc}}$
 and also trivially $m^{\alpha} \rightarrow m$.
The left-hand side of (3.1) converges to
$\text{Re} ((\text{sgn}\, u)[H_{A,m} -m]u)$
in $L^1_{\text{\rm loc}}$, while the right-hand side converges
to $[H_{0,m} - m]|u|$ in distributional sense, so that
 we have shown the desired inequality (1.3).

\medskip
(iii) {\it The case where $m=0$ and $\alpha=1$}. This follows from
the case (ii) for $m>0$, i.e. by letting $m \downarrow 0$ in the equality
(1.3) with $m>0$.
To see this, let $u \in L^2({\Bbb R}^d)$ with
$H_{A,0} u \in L^1_{\text{\rm loc}}({\Bbb R}^d)$.
Then, noting that $H_{A,0} = |-i\nabla-A|$, we see 
by the argument done around (2.1)
that the domains of the operators $H_{A,m}$ and $H_{A,0}$ coincide.
We also see that $H_{A,0} u \in L^1_{\text{\rm loc}}({\Bbb R}^d)$ with
$u \in L^2({\Bbb R}^d)$ implies $H_{A,m} u \in L^1_{\text{\rm loc}}({\Bbb R}^d)$.
In fact, we can show the following fact.

\begin{lemma}
Let $u \in L^2({\Bbb R}^d)$. Then
$H_{A,m} u \in L^1_{\text{\rm loc}}({\Bbb R}^d)$ if and only if
$H_{A,0} u \in L^1_{\text{\rm loc}}({\Bbb R}^d)$.
In fact, for $\psi \in C_0^{\infty}({\Bbb R}^d)$ it holds that
\begin{equation}
\big|\|\psi H_{A,m}u \|_{L^1} - \|\psi H_{A,0}u\|_{L^1}\big|
\leq C(d) m^2\|\psi\|_{L^{\frac{2d}{d+2}}}\|u\|_{L^2}
\end{equation}
with a constant $C(d)$ depending only on $d$.
\end{lemma}

\bigskip
{\it Proof}. We have for $\phi \in C_0^{\infty}({\Bbb R}^d)$
\begin{align}
 H_{A,m}\phi - H_{A,0}\phi
&= ((\!-\nabla -A)^2 +m^2)^{1/2}\phi\, -\, |\!-i\nabla -A|\phi \nonumber\\
&= \big[((\!-\nabla -A)^2 + \theta m^2)^{1/2}\phi\big]_{\theta=0}^{\theta=1}
 = \int_0^1 \frac{d}{d\theta}
           \big[((\!-\nabla -A)^2 + \theta m^2)^{1/2}\phi\big]\, d\theta
                                                               \nonumber\\
&= \frac{m^2}{2}
    \int_0^1 ((-\nabla -A)^2 + \theta m^2)^{-1/2}\phi\, d\theta.
\end{align}
Multiply (3.26) by $\psi \in C_0^{\infty}({\Bbb R}^d)$ with $\psi(x) \geq 0$,
and integrate the absolute value in $x$, then we have
\begin{align}
 \big\|\psi H_{A,m}\phi - \psi H_{A,0}\phi\big\|_{L^1}
&\leq \frac{m^2}{2}
    \int_0^1 \big\|\psi ((-\nabla -A)^2 + \theta m^2)^{-1/2}\phi\big\|_{L^1}\,
                                                        d\theta \nonumber\\
&\leq \frac{m^2}{2}
    \int_0^1 \big\|\psi(-\Delta + \theta m^2)^{-1/2}|\phi|\big\|_{L^1}\,
                                                        d\theta \nonumber\\
&=  \frac{m^2}{2}
    \int_0^1 \int_{{\Bbb R}^d}
    \big[\psi(-\Delta + \theta m^2)^{-1/2}|\phi|\big](x)\,dxd\theta\,,
\end{align}
where the second inequality is due to Lemma 2.1 (i) with
$\beta= \frac12$ and $p=1$.
Note also that the operator $(-\Delta + m^2)^{-1/2}$ in (3.27) has the
following positive integral kernel:
\begin{equation}
 (-\Delta + m^2)^{-1/2}(x)
= \frac{2 m^{d-1}}{(2\pi)^{\frac{d+1}{2}}}
    \frac{K_{(d-1)/2}(m|x|)}{(m|x|)^{(d-1)/2}}, \qquad  m>0,
\end{equation}
with $K_{\nu}(\tau)$ the modified Bessel function of the third kind
of order $\nu$, which was also referred to around (2.8)/(2.9).
In fact, using the expression (2.9) for the integral kernel of
$e^{-tH_{0,m}} = e^{-[-\Delta + m^2]^{-1/2}}$ and integrating it
in $t$ on $(0, \infty)$, we have
\begin{align*}
 (-\Delta + m^2)^{-1/2}(x)
&= \int_0^{\infty} k_0^{m,1}(t,x)\cdot e^{-mt}\,dt\\
&= \int_0^{\infty} 2\big(\frac{m}{2\pi}\big)^{\frac{d+1}{2}}
 \frac{tK_{(d+1)/2}(m(x^2+t^2)^{1/2})}{(x^2+t^2)^{(d+1)/4}}\, dt.
\end{align*}
Change the variables $\tau = m(x^2+ t^2)^{1/2}$, so that
$2tdt = \frac{2\tau}{m^2}d\tau$, and use
$\frac{d}{\tau d\tau}\frac{K_{\nu}(\tau)}{\tau^{\nu}}
            = -\frac{K_{\nu+1}(\tau)}{\tau^{\nu+1}}$,
then we see the last member above be equal to
\begin{align*}
\int_{m|x|}^{\infty} \frac{m^{\frac{d+1}{2}}}{(2\pi)^{\frac{d+1}{2}}}
      \frac{K_{(d+1)/2}(\tau)}{(\tau/m)^{\frac{d+1}{2}}}
      \frac{2\tau}{m^2}d\tau
&= -\frac1{(2\pi)^{\frac{d+1}{2}}}\int_{m|x|}^{\infty}  m^{d+1}
   \frac{d}{\tau d\tau}\Big[\frac{K_{(d-1)/2}(\tau)}{\tau^{\frac{d-1}{2}}}\Big]
      \frac{2\tau}{m^2}d\tau\\
&= \frac{2m^{d-1}}{(2\pi)^{\frac{d+1}{2}}}
      \frac{K_{(d-1)/2}(m|x|)}{(m|x|)^{\frac{d-1}{2}}},
\end{align*}
which yields (3.28).

Since it holds that
$0< K_{\nu}(\tau) \leq C
[\tau^{-\nu}\vee \tau^{-1/2}]e^{-\tau},\,\, \tau>0$
with a constant $C>0$ when $\nu >0$, we obtain
$$
 \frac{K_{(d-1)/2}(\theta^{1/2}m|x|)}{(\theta^{1/2}m|x|)^{\frac{d-1}{2}}}
 \leq C \frac{1}{(\theta^{1/2}m|x|)^{(d-1)}}.
$$
Then we see from (3.27) by the Hardy--Littlewood--Sobolev inequality
(e.g. [LLos01, Chap.4, Sect. 4.3]),
noting $p= \frac{2d}{d+2}$ satisfies the relation
$\frac1{p} + \frac{d-1}{d} + \frac12 = 2$,
\begin{align}
&\big\|\psi H_{A,m}\phi - \psi H_{A,0}\phi \big\|_{L^1} \nonumber\\
&\leq \frac{m^{d+1}}{(2\pi)^{\frac{d+1}{2}}}
   \int_0^1 d\theta\, \theta^{(d+1)/2}
   \int\int_{{\Bbb R}^d\times{\Bbb R}^d} \psi(x)
   \frac{K_{(d-1)/2}(\theta^{1/2}m|x-y|)}{(\theta^{1/2}m|x-y|)^{\frac{d-1}{2}}}
    |\phi(y)|\,dxdy                                     \nonumber\\
&\leq\frac{C(d)}2\frac{m^{d+1}}{(2\pi)^{\frac{d+1}{2}}}
   \int_0^1 d\theta\, \theta\, m^{-(d-1)}
      \|\psi\|_{L^{\frac{2d}{d+2}}}\|\phi\|_{L^2}       \nonumber\\
&= C(d) \frac{m^2}{(2\pi)^{\frac{d+1}{2}}}
  \|\psi\|_{L^{\frac{2d}{d+2}}}\|\phi\|_{L^2}
\end{align}
with a constant $C(d)>0$ depending on $d$.

Now, to show the desired inequality (3.25),
 let $u \in  L^2({\Bbb R}^d)$ and assume that either
$H_{A,m} u$ or $H_{A,0} u$ in $L^1_{\text{\rm loc}}({\Bbb R}^d)$,
consider, for instance, the latter $H_{A,0} u \in L^1_{\text{\rm loc}}({\Bbb R}^d)$.
There exists a sequence $\{\phi_n\}_{n=1}^{\infty}$ in
$C_0^{\infty}({\Bbb R}^d)$ convergent to $u$ in $L^2$ as $n\to \infty$.
We see by (3.29) that
$\{ (\psi H_{A,m} -\psi H_{A,0})\phi_n \}_{n=1}^{\infty}$
is a Cauchy sequence in $L^1$, so that
there exists $v \in L^1({\Bbb R}^d)$ to which it converges
in $L^1$:
$$
 (\psi H_{A,m} -\psi H_{A,0})\phi_n \rightarrow v,
 \qquad n\rightarrow \infty.
$$
{ Since $\psi D[H_{A,0}] \subseteq D[H_{A,0}]$, we see
$\{\psi H_{A,0}\phi_n \}$ converge to
$\psi H_{A,0}u \in L^1({\Bbb R}^d)$
in the weak topology defined by the dual pairing
$\langle L^1({\Bbb R}^d), D[H_{A,0}] \rangle$.
So $\{\psi H_{A,m}\phi_n \}$ becomes a Cauchy sequence
also in this weak topology $\sigma (L^1({\Bbb R}^d), D[H_{A,0}])$,
converging to $v - \psi H_{A,0}u$, which also belongs to
$L^1({\Bbb R}^d)$.
Therefore the existing limit of
$\{\psi H_{A,m}\phi_n \}$ should be written as $\psi H_{A,m}u$
to satisfy
$$
  v = \psi H_{A,m}u + \psi H_{A,0}u.
$$
}
Thus we have seen (3.29) implies
\begin{equation}
\|\psi H_{A,m}u - \psi H_{A,0}u\|_{L^1}
\leq C(d)m^2 \|\psi\|_{L^{\frac{2d}{d+2}}}\|u\|_{L^2}.
\end{equation} 
Hence we have by the triangle inequality
$\big| |a| - |b| \big| \leq \big|a-b\big|$
we have (3.25).
This shows (3.25) for the general $u$, ending the proof of
Lemma 3.3. \qed

\bigskip
Finally, we come back to the proof of Theorem 1.1 to continue
the case (iii) {\it The case where $m=0$ and $\alpha=1$}.
We show that, as $m \downarrow 0$, the left-hand side and the right-hand side
of (1.3) with $m>0$ converge to those with $m=0$.

As for the left-hand side, the sequence $\{\,\|[H_{A,m}-m]u\|_{L^2}^2\,\}$ of 
quadratic forms is increasing as $m$ decreases and converges to $\|H_{A,0} u\|_{L^2}^2$
as $m \downarrow 0$, because
$$
[H_{A,m}-m] = \frac{(-i\nabla-A)^2}{H_{A,m} +m}
          \leq \frac{(-i\nabla-A)^2}{H_A^{m'} +m'} = [H_A^{m'}-m']
          \leq H_{A,0} =|-i\nabla-A|
$$
for $m \geq m' >0$.
This shows the  convergence of the left-hand side of (1.3).
As for the right-hand side, it is easy to see that, as $m \downarrow 0$,
$H_{0,m}|u| \equiv (-\Delta+m^2)^{\frac12}|u|$ converges to
$H_0^0|u| \equiv (-\Delta)^{\frac12}|u|$ in the distribution sense, because
one can show that, for any $\psi \in C_0^{\infty}({\Bbb R}^n)$,
$\{H_{0,m} \psi\}$ converges to $H_0^0 \psi$  as $m \downarrow 0$,
by using
their integral operator representation formula (2.6) with $\alpha=1$;
in fact,
it is due to the convergence of the L\'evy measure $n^{m,1}(dy)$
to the L\'evy measure $n^{0,1}(dy)$ on ${\Bbb R}^d \setminus \{0\}$,
which amounts to the same thing as, observing (2.9),
the convergence of density $n^{m,1}(y)$ to density $n^{0,1}(y)$.
This shows the case $m=0$, completing the proof of Theorem 1.1.
\qed

\bigskip\noindent
{\it Remark}.
From the proof of Theorem 1.1 above, in  particular, the step (i-II), which
relies on Lemma 3.1, we see Theorem 1.1 (Kato's inequality) also hold for
$(H_{A,m})^{\alpha},\,\,(H_{0,m})^{\alpha}$ in place of
$H_{A,m},\,\,H_{0,m}$ with $0<\alpha<1$, that is, (3.1) hold
for $0<\alpha<1$ if $u \in L^2({\Bbb R}^d)$ with
$(H_{A,m})^{\alpha}u \in [L^1_{\text{\rm loc}}({\Bbb R}^d)]^d$.
As a result, Theorem 1.2 (Diamagnetic inequality) also holds for
$(H_{A,m})^{\alpha},\,\,(H_{0,m})^{\alpha}$.

\subsection{Proof of Theorem 1.2}
This has already been implicitly shown in the proof of Theorem 1.1.
In fact, by the same argument used to get (3.20) from (3.18), (3.19),
{\it even for all} $0< \alpha \leq 1$, we have for
$f,\, g \in C_0^{\infty}({\Bbb R}^d)$
\begin{align*}%
 |(f, e^{-t[(H_{A,m})^{\alpha}-m^{\alpha}]}g)|
&\leq e^{m^{\alpha}t}\int_0^{\infty}
             f_{t,\frac{\alpha}2}(\lambda)|(f,e^{-\lambda S}g)|d\lambda\\
&\leq e^{m^{\alpha}t}\int_0^{\infty}
             f_{t,\frac{\alpha}2}(\lambda)(|f|,|e^{-\lambda S}g|)d\lambda\\
&\leq e^{m^{\alpha}t}\int_0^{\infty}
             f_{t,\frac{\alpha}2}(\lambda)(|f|,e^{-\lambda (-\Delta+m^2)}
                |g|)d\lambda\\
&= (|f|, e^{-t[(H_{0,m})^{\alpha}-m^{\alpha}]}|g|).
\end{align*}%
Then this is of course also valid for $f,\, g \in L^2({\Bbb R}^d)$.
\qed.

\section{Concluding Remarks} 
In the literature there are three kinds of relativistic Schr\"odinger operators
for a spinless particle of mass $m\geq 0$ corresponding to
the classical relativistic Hamiltonian symbol
$
\sqrt{(\xi-A(x))^2 +m^2}
$
with {magnetic} vector potential $A(x)$, dependent on {\sl how to quantize
this symbol}. One of them is of course the one $H_{A,m}$ in (1.1)
treated in this paper, and the other two are defined
as pseudo-differential operators, differing from $H_{A,m}$ defined
as an operator-theoretical square root.
In [I12, I13], their common and different properties were
discussed mainly in connection with the
corresponding path integral representations for their semigroups.

The other two relativistic Schr\"odinger operators
are defined by oscillatory integrals, for
$f \in C_0^{\infty}({\Bbb  R}^d)$,  as
\begin{align}%
(H^{(1)}_{A,m} f)(x)
:&\!= \tfrac1{(2\pi)^{d}} \int\!\!\int_{{\Bbb  R}^d\times {\Bbb  R}^d}
   e^{i(x-y)\cdot\xi}
 \sqrt{\Big(\xi-A\big(\tfrac{x+y}{2}\big)\Big)^2 +m^2}\,
   f(y) dyd\xi \nonumber\\
&\!= \tfrac1{(2\pi)^{d}} \int\!\!\int_{{\Bbb  R}^d\times {\Bbb  R}^d}\,
   e^{i(x-y)\cdot (\xi+A(\tfrac{x+y}{2}))}
  \sqrt{\xi^2 +m^2} f(y) dyd\xi \,;  \\
(H^{(2)}_{A,m}f)(x)
:&\!= \tfrac1{(2\pi)^{d}} \int\!\!\int_{{\Bbb  R}^d\times {\Bbb  R}^d}
   \!\!e^{i(x-y)\cdot\xi}
 \sqrt{\big(\xi-
   \int_0^1 A((1-\theta)x+\theta y)d\theta \big)^2 +m^2}\,
   f(y) dyd\xi \nonumber\\
&\!= \tfrac1{(2\pi)^{d}}\! \int\!\!\int_{{\Bbb  R}^d\times {\Bbb  R}^d}
   \!\!e^{i(x-y)\cdot (\xi +\int_0^1 A((1-\theta)x+\theta y)d\theta)}
   \sqrt{\xi^2 +m^2}f(y) dyd\xi. 
\end{align}%
Equality (4.1) is a {\it Weyl pseudo-differential operator} with mid-point prescription
given in [ITa86] (also [I89], [NaU90]) and (4.2)
 a modification of (4.1) given in [IfMP07].
Note that these  two $H^{(1)}_{A,m}$ and $H^{(2)}_{A,m}$ are denoted
in [I12, I13] by slightly different notations $H^{(1)}_{A}$
and $H^{(2)}_{A}$, respectively, while our $H_{A,m}$ in (1.1) by $H^{(3)}_{A}$.

What in this section we should like to call attention to is
that Kato's inequality of distributional form
was {\it missing} for $H^{(3)}_{A,m}$ or our
$H_{A,m}$ in (1.1), although there already exist for the other two
$H^{(1)}_{A,m}$ in (4.1), $H^{(2)}_{A,m}$ in (4.2), indeed it was
shown for $H^{(1)}_{A,m}$ in [I89, ITs92] under some suitable conditions
on $A(x)$  (which differ from $A \in L^2_{\text{\rm loc}}$),
and to be shown in the same way for $H^{(2)}_{A,m}$ (cf. [I13]).
Therefore, at least
the case of Theorem 1.1 with $A=0$ turns out to have already been known.

\medskip
Let us briefly mention here some known facts for
$H^{(1)}_{A,m}$, $H^{(2)}_{A,m}$ and $H^{(3)}_{A,m}$.

1$^{\circ}$. With reasonable conditions on $A(x)$, they
all define {\it selfadjoint} operators  in $L^2({\Bbb  R}^d)$,
which are bounded below.
For instance, they become selfadjoint operators defined as quadratic forms,
for $H^{(1)}_{A,m}$ and $H^{(2)}_{A,m}$, when
$A \in L_{\text{\rm loc}}^{1+\delta}({\Bbb  R}^d; {\Bbb  R}^d)$
for some $\delta>0$
(cf. [I89, I13], [IfMP07]), while for $H^{(3)}_{A,m}$,
when $A \in L_{\text{\rm loc}}^2({\Bbb  R}^d; {\Bbb  R}^d)$
(e.g. [CFKiSi87, pp.8--10] or [I13]).

Furthermore, they are bounded below by the {\sl same lower bound},
in particular,
$$
 H^{(j)}_{A,m} \geq m\,, \qquad j= 1, \,2, \,3.
$$

2$^{\circ}$. $H^{(2)}_{A,m}$ and $H^{(3)}_{A,m}$ are covariant under
gauge transformation,
i.e. for every
$\varphi \in {\cal S}({\Bbb  R}^d)$ it holds that
$\,\,H^{(j)}_{A+\nabla \varphi} = e^{i\varphi}H^{(j)}_{A,m}e^{-i\varphi}$,
$j=2,3$. However, $H^{(1)}_{A,m}$ is not.

3$^{\circ}$. All these three operators are different in general,
but coincide, if $A(x)$ is
linear in $x$, i.e. if $A(x) = \dot{A}\cdot x\,$ with $\dot{A}: d\times d$
real symmetric {\sl constant matrix}, then
$H^{(1)}_{A,m} = H^{(2)}_{A,m} = H^{(3)}_{A,m}$.
So, this holds for uniform magnetic fields with $d=3$.

\bigskip\bigskip\noindent
{\bf{\Large Appendix A}} 

\bigskip\noindent
Our aim is to derive the following expressions for integral kernel
$k_0^{m,\alpha}(t,x)$ of semigroup $e^{-t[(H_{0,m})^{\alpha}-m^{\alpha}]}$ and
 density (function) $n^{m,\alpha}$ of L\'evy measure $n^{m,\alpha}(dy)$
for $0<\alpha \leq 1$, which are mentioned around formulas 
(2.7), (2.8)/(2.9):
\begin{align}%
k_0^{m,\alpha}(t,x)
&= \frac{e^{m^{\alpha}t}}{\pi(2\pi)^{\frac{d}2}|x|^{\frac{d}2-1}}
   \int_0^{\infty} e^{-tr^{\frac{\alpha}2}\cos \tfrac{\alpha}2\pi}
 \sin(tr^{\frac{\alpha}2}\sin \tfrac{\alpha}2\pi)(m^2+r)^{\frac12(\frac{d}2-1)}
                                                                 \nonumber\\
   &\qquad\qquad\qquad\times K_{\frac{d}2-1}((m^2+r)^{\frac12}|x|)\,dr \tag{A.1}\\
n^{m,\alpha}
&= \frac{2^{1+\frac{\alpha}2}\sin \big(\frac{\alpha}2\pi\big)
   (2\pi)^{\frac{\alpha}2}\Gamma(\tfrac{\alpha}2+1)}{\pi}
   \big(\frac{m}{2\pi}\big)^{\frac{d+\alpha}2}
   \frac{K_{\frac{d+\alpha}2}(m|x|)}{|x|^{\frac{d+\alpha}2}} \tag{A.2}.
\end{align}%
Equality (A.2) is essentially the same as $\nu^m$
in [ByMaRy09, (2.7), p.4877], which is established for the heat semigroup
$e^{-t[(-\Delta+m^{\alpha/2})^{\alpha/2}-m]}$ instead of our
$e^{-t[(H_{0,m})^{\alpha}-m^{\alpha}]}$. Indeed, putting in (A.2)
$m= {m'}^{\frac{1}{\alpha}}$ to rewrite it with Euler's reflection formula
 $\Gamma(z)\Gamma(1-z) = \frac{\pi}{\sin(\pi z)}$ yields
eq. (2.7) in this reference with $m$ replaced by $m'$.

To show (A.1) and (A.2), we use another formula (3.18)/(3.19)
to express the semigroup
$e^{-t(H_{0,m})^{\alpha}}
                     \equiv e^{-t(-\Delta+m^2)^{\frac{\alpha}2}}\,\,
(0<\alpha \leq 1)$ for the fractional power:
\begin{align*}%
(e^{-t(H_{0,m})^{\alpha}}u)(x)
&= \int_{{\Bbb R}^d}\Big(\int_0^{\infty} f_{t,\frac{\alpha}2}(s)
e^{-s(-\Delta+m^2)}ds\,u\Big)(y) dy,\\
\intertext{where $e^{-t(-\Delta+m^2)}$ is the heat semigroup
multiplied by $e^{-m^2t}$:}%
(e^{-t(-\Delta+m^2)}u)(x)
&= \frac1{(4\pi t)^{d/2}}\int_{{\Bbb R}^d} e^{-m^2t- \frac{(x-y)^2}{4t}}u(y)dy.
\end{align*}%
The $f_{t,\frac{\alpha}2}(s)$ in (3.19) is rewritten as%
$$%
f_{t,\frac{\alpha}2}(s)
= \frac1{\pi}\int_0^{\infty} e^{sr\cos\theta-tr^{\alpha}\cos \frac{\alpha}2\theta}
                   \sin(sr\sin\theta-tr^{\frac{\alpha}2}
                   \sin \tfrac{\alpha}2\theta +\theta)dr
 \quad (t>0, \,\,s\geq 0),
$$%
where the integration path is deformed to the union of two paths
$re^{-i\theta} (-\infty <-r<0)$ and $re^{i\theta} (0<r< \infty)$, where
$\frac{\pi}2 \leq \theta \leq \pi$
(see [Y78, IX, 11, pp.259--263]).

Then we take $\theta=\pi$ to have
\begin{align*}%
&(e^{-t(H_{0,m})^{\alpha}})(x)\\
&= \int_0^{\infty} ds \frac1{\pi}
   \frac1{(4\pi s)^{\frac{d}2}} e^{-m^2s- \frac{x^2}{4s}}
   \int_0^{\infty} e^{-sr-tr^{\frac{\alpha}2}\cos \tfrac{\alpha}2\pi}
                   \sin(tr^{\frac{\alpha}2}\sin \tfrac{\alpha}2\pi)dr\\
&= \frac1{\pi(4\pi)^{\frac{d}2}}
   \int_0^{\infty}dr\,  e^{-tr^{\frac{\alpha}2}\cos \tfrac{\alpha}2\pi}
     \sin(tr^{\frac{\alpha}2}\sin \tfrac{\alpha}2\pi)(m^2+r)^{\frac{d}2-1}
 \int_0^{\infty}\frac{e^{-s- \frac{(m^2+r)x^2}{4s}}}{s^{\frac{d}2}}ds \\
&= \frac{1}{\pi(2\pi)^{\frac{d}2}|x|^{\frac{d}2-1}}
  \int_0^{\infty} e^{-tr^{\frac{\alpha}2}\cos \tfrac{\alpha}2\pi}
 \sin(tr^{\frac{\alpha}2}\sin \tfrac{\alpha}2\pi)(m^2+r)^{\frac12(\frac{d}2-1)}\\
  &\qquad\qquad\qquad\qquad\qquad\qquad\qquad\qquad\qquad
   \times K_{\frac{d}2-1}((m^2+r)^{\frac12}|x|)\,dr,
\end{align*}%
where we have used the representation formula of the modified Bessel function
of the third kind, $K_{\nu}(z)$ [GrR94, sect. 8.432. 6, p.969]:
$$
K_{\nu}(z) = \tfrac12\big(\tfrac{z}{2}\big)^{\nu}\int_0^{\infty}
           e^{-t-\frac{z^2}{4t}}\, t^{-\nu-1} dt,
\quad \nu > -\tfrac12,\,\, z>0.
$$
It follows that the integral kernel $k_0^{m,\alpha}(t,x)$ of the semigroup
$e^{-t[(H_{0,m})^{\alpha}-m^{\alpha}]}$ turns out
\begin{align*}%
k_0^{m,\alpha}(t,x)
:&= e^{-t[(H_{0,m})^{\alpha}-m^{\alpha}]}(x)\\
&= \frac{e^{m^{\alpha}t}}{\pi(2\pi)^{\frac{d}2}|x|^{\frac{d}2-1}}
   \int_0^{\infty} e^{-tr^{\frac{\alpha}2}\cos \tfrac{\alpha}2\pi}
 \sin(tr^{\frac{\alpha}2}\sin \tfrac{\alpha}2\pi)(m^2+r)^{\frac12(\frac{d}2-1)}\\
   &\qquad\qquad\qquad\times K_{\frac{d}2-1}((m^2+r)^{\frac12}|x|)\,dr\,.
\end{align*}%
This shows (A.1).

Next, we have
\begin{eqnarray*}%
&&\frac{d}{dt}k_0^{m,\alpha}(t,x)\\
&=&\frac{1}{\pi(2\pi)^{\frac{d}2}|x|^{\frac{d-1}2}}\int_0^{\infty}dr\,
     \frac{d}{dt}\Big[e^{t(m^{\alpha}-r^{\frac{\alpha}2}\cos \tfrac{\alpha}2\pi)}
\sin (tr^{\frac{\alpha}2}\sin \tfrac{\alpha}2\pi) \\
  &&\qquad\qquad\qquad\qquad\qquad\qquad
  \times (m^2+r)^{\frac12(\frac{d}2-1)}K_{\frac{d}2-1}((m^2+r)^{\frac12}|x|)\Big]\\
&=& \frac{1}{\pi(2\pi)^{\frac{d}2}|x|^{\frac{d-1}2}}\int_0^{\infty}dr\,
 \Big[(m^{\alpha}-r^{\frac{\alpha}2}\cos \tfrac{\alpha}2\pi)
    \sin (tr^{\frac{\alpha}2}\sin \tfrac{\alpha}2\pi)
 + r^{\frac{\alpha}2}\sin \tfrac{\alpha}2\pi
       \cos (tr^{\frac{\alpha}2}\sin \tfrac{\alpha}2\pi)\Big]\\
 &&\qquad\qquad\qquad\qquad\qquad\qquad
  \times  e^{t(m^{\alpha}-r^{\frac{\alpha}2}
       \cos \frac{\alpha}2\pi)}(m^2+r)^{\frac12(\frac{d}2-1)}
   K_{\frac{d}2-1}((m^2+r)^{\frac12}|x|).\\
\end{eqnarray*}
Then by the fact (2.7), we have, as $t \downarrow 0$,
\begin{eqnarray*}%
n^{m,\alpha}(t,x) &= \frac1{t}k_0^{m,\alpha}(t,x)\\
&& \rightarrow
\frac{d}{dt}k_0^{m,\alpha}(t,x)\Big|_{t=0}=: n^{m,\alpha}(x) \\
 &\quad&\, = \frac{\sin \frac{\alpha}2\pi}{\pi(2\pi)^{\frac{d}2}{|x|^{\frac{d}2-1}}}
 \int_0^{\infty}dr\, (m^2+r)^{\frac12(\frac{d}2-1)}r^{\frac{\alpha}2}
   K_{\frac{d}2-1}((m^2+r)^{\frac12}|x|).
\end{eqnarray*}%
Here the integral on the last member above is equal to
\begin{eqnarray*}%
&&\,\, \int_0^{\infty}\, (m^2+\tau^2)^{\frac12(\frac{d}2-1)}\tau^{\alpha}
   K_{\frac{d}2-1}((m^2+\tau^2)^{\frac12}|x|)\, 2\tau d\tau
   \quad(r := \tau^2)\\
&=& 2\int_m^{\infty}\, a^{\frac{d}2-1} (a^2-m^2)^{\frac{1+\alpha}2}
   K_{\frac{d}2-1}(a|x|)\, \frac{a}{(a^2-m^2)^{\frac12}}da
   \quad(a := (m^2+\tau^2)^{\frac12}\\
&=& \frac2{|x|^{\frac12}}\int_m^{\infty}\, a^{\frac{d-1}2} (a^2-m^2)^{\frac{\alpha}2}
   K_{\frac{d}2-1}(a|x|)\, (a|x|)^{\frac12}da. \\
\end{eqnarray*}%
Then we use the following formula
[EMOT54, Chap.X ({\it K-Transforms}), 10.2.(13), p.129]:
\begin{eqnarray*}%
\int_a^{\infty} x^{\frac12-\nu}(x^2-a^2)^{\mu}K_{\nu}(xy)(xy)^{\frac12}dx
&=& 2^{\mu}a^{\mu-\nu+1}y^{-\mu-\frac12}\Gamma(\mu+1)
   K_{\mu-\nu+1}(ay),\\
&&\qquad y >0, \,\, \mu >-1,
\end{eqnarray*}%
reading with $\mu=\frac{\alpha}2, \, -\nu = \frac{d}2-1$ and
with ``$\nu$" in place of ``$-\nu$"because $K_{-\nu}(\tau) = K_{\nu}(\tau)$.
to finally obtain
\begin{eqnarray*}%
n^{m,\alpha}(x)
&=& \frac{\sin \frac{\alpha}2 \pi}{\pi(2\pi)^{\frac{d}2}|x|^{\frac{d}2-1}}
   2^{\frac{\alpha}2 +1}m^{\frac{d+\alpha}2}
   |x|^{-(\frac{\alpha}2+1)}\Gamma(\tfrac{\alpha}2+1)K_{\frac{d+\alpha}2}(m|x|)\\
&=& \frac{2^{1+\frac{\alpha}2}\sin \big(\frac{\alpha}2\pi\big)
   (2\pi)^{\frac{\alpha}2}\Gamma(\tfrac{\alpha}2+1)}{\pi}
   \big(\frac{m}{2\pi}\big)^{\frac{d+\alpha}2}
   \tfrac{K_{\frac{d+\alpha}2}(m|x|)}{|x|^{\frac{d+\alpha}2}}.
\end{eqnarray*}%
If $\alpha=1$, this expression reduces to
$$
n^{m,1}(x)
= 2\big(\frac{m}{2\pi}\big)^{\frac{d+1}2}
    \tfrac{K_{\frac{d+1}2}(m|x|)}{|x|^{\frac{d+1}2}},
$$
which is nothing but the first formula of (2.9),
and we see that $n^{m,\alpha}(x)$ tends to $n^{m,1}(x)$,
as $\alpha \uparrow 1$ since
$\Gamma(\frac12+1) = \frac{\pi^{\frac12}}2$.

\bigskip\bigskip\bigskip\noindent
{\bf Acknowledgments}

\medskip\noindent
The work by FH is supported in part by JSPS Grant-in-Aid for Challenging
Exploratory Research 15K13445, and
the work by TI supported in part by JSPS Grant-in-Aid for Scientific
Research 23540191 and 25400166.

The authors would like to thank an anonymous referee for very careful
reading of the manuscript and valuable comments improving the paper.


\end{document}